\documentclass[11pt]{article}

\newtheorem{rem}{Remark}[section]
\newtheorem{teo}{Theorem}[section]
\newtheorem{defin}{Definition}[section]
\newtheorem{prop}{Proposition}[section]
\newtheorem{coro}{Corollary}[section]
\newtheorem{lemma}{Lemma}[section]
\newtheorem{es}{Example}[section]
\def\proof{{\bf Proof:}\ }
\def\endproof{{\mbox{}\nolinebreak\hfill\rule{2mm}{2mm}\par\medbreak} }

\def\eq#1{(\ref{#1})}
\def\neweq#1{\begin{equation}\label{#1}}
\def\endeq{\end{equation}}

\newtheorem{remark}{Remark}[section]
\newtheorem{problem}{Problem}[section]

\font\tenmsa=msam10 scaled 1200 \font\sevenmsa=msam7 scaled 1200
\font\fivemsa=msam5 scaled 1200 \font\tenmsb=msbm10 scaled 1200
\font\sevenmsb=msbm7 scaled 1200 \font\fivemsb=msbm5 scaled 1200
\newfam\msafam
\newfam\msbfam
\textfont\msafam=\tenmsa \scriptfont\msafam=\sevenmsa
\scriptscriptfont \msafam=\fivemsa \textfont\msbfam=\tenmsb
\scriptfont\msbfam=\sevenmsb \scriptscriptfont\msbfam=\fivemsb

\newcommand{\amsb}  {\fam\msbfam}

\def\bbbr{{\amsb R}}
\def\bbbn{{\amsb N}}

\def\supp(#1){\{#1>0\}}

\def\a{\alpha}

\def\eps{\varepsilon}
\def\l{\lambda}
\def\O{\Omega}

\def\G{\Gamma}

\def\R{\bbbr}

\def\f{\varphi}

\def\eps{\varepsilon}
\def\f{{\varphi}}
\def\div{{\rm div}}
\def\cal#1{{\mathcal #1}}

\def\l{\lambda}
\def\O{\Omega}

\def\f{\varphi}

\def\div{{\rm div}}

\title{\sc A variational problem for the spatial segregation of reaction--diffusion systems
\footnote{\sc Work partially supported by
MIUR, Project ``Metodi Variazionali ed Equazioni Differenziali
Non Lineari'' }}
\author{Monica Conti
\footnote{%
Dipartimento di Matematica  del Politecnico di Milano,
Piazza Leonardo da Vinci, 32;  20133 Milano, Italy.
email: \emph{monica@mate.polimi.it}}~,
~Susanna Terracini%
\footnote{%
Dipartimento di Matematica e Applicazioni
Università degli Studi di Milano-Bicocca
Via Bicocca degli Arcimboldi, 8; 20126 Milano, Italy.
email: \emph{suster@matapp.unimib.it}
}
~and Gianmaria Verzini
\footnote{%
Dipartimento di Matematica  del Politecnico di Milano,
Piazza Leonardo da Vinci, 32;  20133 Milano, Italy.
email: \emph{giaver@mate.polimi.it}}
}

\date{\today}

\begin{document}
\maketitle

\begin{abstract} In this paper we study a class of stationary states
for reaction--diffusion systems of $k\geq 3$ densities having disjoint supports.
For a class of segregation states governed by
a variational principle we prove existence and provide conditions for uniqueness.
Some qualitative properties and the local regularity  both of the densities and of
their free boundaries are established in
the more general context of a functional class characterized by differential inequalities.
\end{abstract}




\section{Introduction}
The occurrence of nontrivial steady states (pattern formation) for
reaction--diffusion systems has been widely studied in the
literature. Of particular interest is the existence of spatially
inhomogeneous solutions for competition models of Lotka--Volterra
type. This study has been carried out mainly in the case
of two competing species, see e.g. \cite{bb,cl,eflg,gl,kl,lm,mhmm,mm,skt}; in recent years also the case of
three competing densities, which is far more complex, has become
object of an extensive investigation \cite{dd,ddh,lmn,ln}. In most cases, the
pattern formation is driven by the presence of different diffusion
rates when the coefficients of intra--specific and inter--specific
competitions are suitable related. A remarkable
limit case of pattern formation yields to the segregation of
competing species, that is, configurations where different
densities have disjoint habitats; see \cite{d1,dd3,dg1,dhmp,ey,t}. Object of the
present paper is to study a class of possible segregation states,
involving an arbitrary number of competing densities, which are
governed by a minimization principle rather than
competition--diffusion.
 Roughly speaking, we are going to deal
with stationary configurations of $k\geq 2$ densities that
interact only through the boundaries of their nodal sets; the
minimization involves the sum of the internal energies, with the
constraint of being segregated states. In other words, the
supports of the densities have to satisfy a suitable optimal
partition problem in $\bbbr^N$. Precisely, let $\O$
be a bounded open subset of $\bbbr^N$ ($N\geq 2$) and let us call {\it segregated state}
a $k$--uple $U=(u_1,\dots,u_k)\in (H^1(\Omega))^k$ where
$$u_i(x)\cdot u_j(x)=0\hspace{1cm}i\neq j,\mbox{ a.e. } x \in\O.$$
We define the internal energy of $U$ as
$$J(U)=\sum_{i=1,\cdots,
k}\left\{\int_\Omega\left(\frac12\,d_i^{\,2}(x) \,|\nabla
u_i(x)|^2- F_i(x,u_i(x))\right)dx\right\}, $$  where the $F_i$'s
satisfy suitable assumptions (see (A1), (A2) below). Our first goal is to minimize
$J$  among a class of segregated states subject to some boundary and positivity conditions;
next we shall develop a regularity and a free boundary theory for minimizers.
In performing the second goal the main tools will come from recent results that
Caffarelli, Jerison and Koenig \cite{cjk}
(see also \cite{acf, acs, cks, cjk1}  and references therein) have obtained in the study of
free boundaries in other contexts. Recently, free boundary problems have been studied in
connection with the asymptotic behaviour of some models of population dynamics with diffusion
as in Dancer and Du (\cite{dd2}) and Dancer, Hilhorst et al in (\cite{dhmp}), in the case of two competing
species. In a forthcoming paper we shall show how our variational problem appears as a limiting
problem for some classes of competition--diffusion problems.

Our first result
establishes existence for this problem; then we discuss the uniqueness
of the solution. Surprisingly enough, the minimizer can be proven to be unique
 for a large class of Lagrangians. A remarkable fact is that
solutions to this variational problem satisfy extremality
conditions in the form of differential inequalities of special type. Precisely they
belong to the functional class
$${\cal S}=\left\{(u_1,\cdots,u_h)\in
(H^1(\Omega))^k:\,\begin{array}{l}
 u_i\geq0,\, u_i\cdot u_j=0\mbox{ if }i\neq j\\
 -\Delta u_i\leq f_i(x,u_i),\,-\Delta \widehat u_i\geq\widehat f(x,\widehat u_i)
\end{array}\right\}
$$
where $\widehat u_i = u_i-\sum_{h\neq i}u_h$ and $\widehat f(x,\widehat u_i)
=\sum_{j} f_j(x,\widehat u_i)\chi_{\supp(u_j)}$.

A further reason of interest in the class $\cal S$ is that
it contains also the asymptotic limits of the solutions of a large class of competition--diffusion systems when
the inter--specific competition terms tend to infinity.  This will be the object of a forthcoming paper; a
link between some variational problems and competing species
systems has been traced by the authors in \cite{ctv,ctv2}.

An important part of the paper is devoted to study the qualitative
properties exhibited by the segregated states belonging
to the class ${\cal S}$. In particular we shall establish the
local lipschitz continuity both of $U$ and its nodal set; to this aim we will take advantage of some monotonicity formulae as in \cite{acf,cjk}.
Then, for the dimension $N=2$, we develop further our investigation.
Our main result is that, near a zero point, $U$ and its null set
exhibit the same qualitative behavior of  harmonic functions and their nodal sets
(\cite{alex1,hw}).
In particular we prove that
the set of double points (i.e. points where two densities meet) is the union of a finite
number of regular arcs meeting at a finite number of multiple points
(i.e. points where more than two densities meet).
We emphasize that, at a multiple point the densities share the angle in equal parts
and moreover an asymptotic expansion for $U$ is available.

The plan of the paper is the following: in Section \ref{ipotesi} we introduce the basic assumptions
and formulate the variational problem;
the existence of a minimizer is proven in Section \ref{exist-depend}; Section \ref{unicita} deals with the uniqueness of the solution; finally
in Section \ref{stoinesse} the extremality conditions are established.
In Sections \ref{sez.classe.esse} and \ref{sez.classe.esse*} we introduce ${\cal S}$
and a wider functional class ${\cal S}^*$; next in Section \ref{regularity} we prove
the local lipschitz continuity in ${\cal S}^*$ and the global regularity in ${\cal S}$.
In Section \ref{N=2} we establish some qualitative properties of the elements of ${\cal S}$ in dimension $N=2$.

\section{Assumptions and notation}\label{ipotesi}

Let $N\geq 2$; let $\Omega\subset{\bbbr}^N$ be a connected, open bounded domain
with regular boundary $\partial\Omega$. Let  $k\geq 2$
be a fixed integer.
Throughout all the paper we will make the
following set of assumptions (for every $i=1,\dots,k$):
\begin{itemize}
\item
$\phi_i\in H^{1/2}(\partial\Omega)$,
$\phi_i\geq0$, $\phi_i\cdot\phi_j=0$ for $i\neq j$, almost everywhere on
$\partial\Omega$; sometimes such a boundary datum  will by called {\it admissible}
\item $d_i\in
W^{2,\infty}(\Omega)$, $d_i>0$ on $\overline\Omega$
\item  $f_i(x,s):\Omega\times \R^+\to\R$ such that:
\begin{itemize}
\item[(A1)] $f_i(x,s)$ is Lipschitz continuous in $s$, uniformly in $x$
 and $f_i(x,0)\equiv0\label{hip1}$
\item[(A2)] there exists  $b_i\in L^\infty(\Omega)$ such
that both
$$
|f_i(x,s)|\leq b_i(x)s\qquad\forall\, x\in\Omega,\;
s\geq\bar s>>1
$$
and
$$
\int_\Omega \left(\,d_i^2(x)|\nabla
w(x)|^2-b_i(x)w^2(x)\right)dx>0\qquad\forall\,w\in
H^1_0(\Omega)\label{hip3}.
$$
\end{itemize}
\end{itemize}
Every $\phi_i$ will be the boundary trace of a non negative
density $u_i\in H^1(\Omega)$. Moreover, associated to each
density, we consider its diffusion coefficient $d_i$ and its internal potential
$F_i(x,s):=\int_0^s f_i(x,u)du$. \\
We are concerned with the
following variational problem.
\begin{problem}\label{Prb}\rm Let
$${\mathcal U}=\left\{(u_1,\cdots,u_k)\in(H^{1}(\Omega))^k:
u_i|_{\partial\Omega}=\phi_i,u_i\geq0\; \forall\,i=1,\cdots,
k;\;u_j\cdot u_i=0, i\neq j\mbox{ a.e. on }\Omega\right\}.$$ Find the
minimum of the functional
\begin{equation}\label{energia}J(U)=\sum_{i=1,\cdots, k}\left\{\int_\Omega\left(\frac12\,d_i^{\,2}(x)
\,|\nabla u_i(x)|^2- F_i(x,u_i(x))\right)dx\right\}
\end{equation}
where $U\in{\mathcal U}$.
\end{problem}

\begin{remark}
\rm We notice that the case with variable diffusions can be
reduced, by a suitable change of the variables $u_i$'s, to the
case when $d_i\equiv1$ for every $i$. Indeed, let $u_i=v_i/d_i$
(this can be done because the $d_i$'s are strictly positive).
After integrating by parts we obtain
\begin{eqnarray*}J(U)&=&\sum_{i=1,\cdots,
k}\left\{\int_\Omega\left(\frac12 \,|\nabla v_i|^2+\frac12\frac{\Delta d_i}{d_i}v_i^2-
F_i\left(x,\frac{v_i}{d_i}\right)\right)dx-\frac12\int_{\partial\Omega}
\phi_i^2 d_i \frac{\partial
d_i}{\partial\nu}ds\right\}\\&=&\sum_{i=1,\cdots,
k}\left\{\int_\Omega\left(\frac12 \,|\nabla v_i|^2-\widetilde
{F_i}(x,v_i)\right)dx-C\right\}=\widetilde J(V)-C.
\end{eqnarray*}
Exploiting this identity, the reader can easily check that the
validity of the assumptions for the $f_i$'s and for the $\phi$'s
implies the same for the new data  $\widetilde {f_i}$'s and $\widetilde {\phi_i}$'s.
Hence, in what follows, we will choose $d_i\equiv1$ for every $i$, unless otherwise
specified (namely in the results of Section \ref{unicita}).
\end{remark}

\begin{remark}\rm
By our definition, the functions $f_i$'s are defined only for non
negative values of $s$ (recall that our densities $u_i$'s are
assumed non negative); thus we can arbitrarily define such
functions on the negative semiaxis. For the sake of convenience,
when $s\leq0$ we will let $f_i(x,s): =-f_i(x,-s)$. This extension
preserves the continuity, thanks to assumption (A1). In the same
way, each $F_i$ is extended as an even function.
\end{remark}
\small {\bf Notation}
In the following, when not needed, we shall omit the dependence on
the variable $x$.
We use the standard notation  $g^+(x)=\max_{x\in\Omega}(g(x),0)$ and
$g^-(x)=\max_{x\in\Omega}(-g(x),0)$.
Given a $k$--uple $(u_1,\dots,u_k)$ we introduce the ``hat'' operation  as
$$\widehat{u}_i:=u_i-\sum_{j\neq i}u_j.$$
Furthermore, with some abuse of notation,
 we shall use a capital letter to identify both a $k$--uple and the sum of its $k$ components
 (e.g. $U=(u_1,\dots,u_k)$ and $U=\sum_{i=1}^k u_i$ ).
With the notation $(u_{i,n})$ we shall denote the $i$--th component of a sequence of $k$--uples $(U_n)$.
The symbol $\chi_A$ will denote the characteristic function of the set $A$.

\normalsize
\section{Existence of the minimum and continuous dependence}\label{exist-depend}
Our first goal is to prove the existence of at least one minimizer
of Problem \ref{Prb}. Next we shall prove continuous dependence of the minimizers (that do not need
to be unique) with respect to the data. This continuity property will be exploited in the analysis of
 the local properties of the solutions, and prescisey when performing the blow--up
argument.

To this aim, we start observing that our assumptions on $f_i$
imply
\begin{equation}\label{controlprim}
|F_i(x,s)|\leq \frac{b_i(x)}{2}s^2+C|s|\quad\forall
x\in\O,\,\forall s\in\bbbr
\end{equation}
for every $i$. On the other hand, by standard eigenvalues theory, assumptions
(A2) implies that the quadratic form there is an
equivalent norm on $H^1_0(\Omega)$, that is, there exists
$\varepsilon>0$ such that
$$
\int_\Omega (|\nabla
w|^2-b_i(x)w^2)dx\geq\varepsilon\int_\Omega|\nabla w|^2dx,
$$
for every $w\in H^1_0(\Omega)$ and for every $i$. As a consequence
we have the following result:

\begin{teo}\label{esistenza} Under the assumptions of Section \ref{ipotesi}
{\rm Problem \ref{Prb}} has at least one solution.
\end{teo}
\begin{proof} applying \eq{controlprim} we easily obtain that
the $H^1$-continuous functional \eq{energia} is coercive;
indeed, for every $u_i\in H^1(\Omega)$, we have $$\int_\Omega
\left(\frac12|\nabla u_i|^2-F_i(u_i)\right)dx \geq\int_\Omega
\frac12\left(|\nabla u_i|^2-b_iu_i^2\right)dx-
\int_{\Omega}Cu_idx$$
$$\geq\frac{\varepsilon}{2}\int_{\Omega}|\nabla u_i|^2-c(\Omega)
$$

 for some constant $c(\Omega)$. Let us take a minimizing
sequence $(u_{i,n})$ in $\cal U$ of disjoint support functions.
This sequence being $H^1$-bounded by the above inequality, there
exists a subsequence weakly convergent to $u_i$ in $H^1$ and, by
compact injection, in the $L^2$--strong topology; taking possibly
a new subsequence, we  infer almost everywhere convergence in
$\Omega$ of every $u_{i,n}$ to $u_i$ and the limit functions
have obviously disjoint supports. The weak lower semicontinuity
ensures that the weak limit is in fact a minimizer.
\end{proof}

\begin{teo}\label{teo3}
Let $U_n=(u_{1,n},\cdots,u_{k,n})$ be a sequence of solutions to {\rm Problem \ref{Prb}},
respectively with admissible data $(\phi_{1,n},\cdots,\phi_{k,n})$, such that
$$\phi_{i,n}\longrightarrow\phi_i\qquad \mbox{in
}H^{1/2}(\partial\Omega)$$ and potentials
$$F_{i,n}(x,s)\longrightarrow F_i(x,s)\qquad \mbox{in
}C^1(\Omega\times \R)$$ for every $i=1,\cdots,k$. Then,
$$u_{i,n}\longrightarrow u_i\qquad \mbox{in }H^{1}(\Omega)$$for
every $i=1,\cdots,k$ and $U=(u_1,\cdots,u_k)$ solves {\rm Problem \ref{Prb}} with data
$(\phi_1,\cdots,\phi_k)$ and potentials $(F_1,\cdots,F_k)$.\end{teo}

\begin{proof}
first, we observe that $(\phi_1,\cdots,\phi_k)$ is an admissible datum, i.e. the $\phi_j$'s are nonnegative and have
disjoint supports by the strong convergence of $(\phi_{1,n},\cdots,\phi_{k,n})$ in
$H^{1/2}(\partial\Omega)$.

We denote by $U^{*}=(u_1^{*},\cdots,u_k^{*})$ a solution to Problem \ref{Prb} with data
$(\phi_1,\cdots,\phi_k)$, and $(F_1,\cdots,$ $F_k)$.
Consider the minimum levels
$$c_n= J(U_n)\quad\mbox{and}\quad c^{*}= J(U^{*}).$$
Observe that the convergence of the boundary traces $\phi_{i,n}$'s and of the $F_{i,n}$'s, ensures a bound on the sequence $c_n$.
The coercivity of $J$ then yields to a bound on the sequence $\|u_{i,n}\|_{H^1(\Omega)}$; therefore we can assume,
up to a subsequence, that
\begin{eqnarray}\nonumber c_n&\longrightarrow& c^0
\\
u_{i,n}&\longrightarrow&u_i\mbox{ weakly in
}H^1(\Omega).\end{eqnarray}
Furthermore, as a
consequence of the compact injection $H^1\hookrightarrow L^2$, we
have $u_i\cdot u_j=0$, whenever $i\neq j$. Moreover, by the weak continuity of the
trace operator, we obtain $$u_i|_{\partial\Omega}=\phi_i.$$
The
lower weak semicontinuity of the norm implies
$$
\frac12\sum_{i=1 ,  \cdots, k}\int_\Omega|\nabla
u_i|^2\leq\liminf_{n \to\infty}\frac12\sum_{i=1\cdots
k}\int_\Omega|\nabla u_{i,n}|^2
$$
and we also have
$$\sum_{i=1,   \cdots, k}\int_\Omega
F_i(x,u_i(x))=\lim_{n\to\infty}\sum_{i=1,   \cdots, k}\int_\Omega
F_{i,n}(x,u_{i,n}(x)).$$

We observe that the level $$c=\sum_{i=1   \cdots
k}\left\{\frac12\int_\Omega|\nabla u_i|^2-\int_\Omega
F_i(u_i)\right\}$$ is not necessary a minimum level but satisfies,
by the discussion above, the inequalities
$$c_0\geq c\geq c^*.$$

We wish to prove that $c_0=c^*$. Suppose, by contradiction, that $c^{*}<c_0.$  Consider the harmonic
extensions (still denoted with the same symbols) on $\Omega$ of the $\phi_{i,n}$'s and of the $\phi_i$'s and introduce $\psi_{i,n}=\phi_{i,n}-\phi_i$.
Then, by construction
\begin{eqnarray} \label{vanish}\psi_{i,n}&\longrightarrow& 0\,\,\mbox{ in
}H^1(\Omega)\\
\nonumber\left.\psi_{i,n}\right|_{\partial
\Omega}&\longrightarrow& 0\,\,\mbox{ in }H^{1/2}(\partial\Omega).
\end{eqnarray}
We define
\begin{eqnarray*}
w_{i,n}&=&\Big(u_i^{*}+\psi_{i,n}\Big)^+\\v_{i,n}&=&
\Big(w_{i,n}-\sum_{j\neq i}w_{j,n}\Big)^+.\end{eqnarray*}
We observe that $w_{i,n}|_{\partial\Omega}=\phi_{i,n}$; moreover,
since $u_i^*\geq 0$ and by \eq{vanish},
\begin{eqnarray}
\nonumber w_{i,n}&\longrightarrow& u_i^{*}\mbox{
in }H^1(\Omega)\label{eq.g}\\v_{i,n}&\longrightarrow&
u_i^{*}\mbox{ in }H^1(\Omega).\label{eq.h}
\end{eqnarray}
Moreover, since $ w_{i,n}\geq 0$, it is immediate to see that $v_{i,n}\cdot v_{j,n}=0$
if $i\neq j$. Hence it follows from the definition of $c_n$ that
$$\sum_{i=1,\cdots,k}\left\{\frac12\int_\Omega|\nabla
v_{i,n}|^2-\int_\Omega F_{i,n}(v_{i,n})\right\}\geq c_n\;;$$
on the other hand \eq{eq.h} implies $$
\sum_{i=1,\cdots,k}\left\{\frac12\int_\Omega|\nabla v_{i,n}|^2-\int_\Omega
F_{i,n}(v_{i,n})\right\}\longrightarrow c^{*}$$
that implies that $c^{*}\geq c_0$.

Finally, from the equality of the minima $c_0=c=c^*$, we also deduce
the strong convergence of the  $u_{i,n}$'s in $H^1(\Omega)$ and the thesis follows.
\end{proof}
\section{Uniqueness}\label{unicita}
In general we can not expect the minimizer of Problem \ref{Prb} to be unique. A simple
counterexample can be constructed in the following way:
\begin{es} Let $k=2$, $\phi_1\equiv\phi_2\equiv0$,
$f_1(x,s)=f_2(x,s)=min(\lambda s,s^{1/3})$, where
$\lambda>\lambda_1(\Omega)$, the first eigenvalue of the Laplace
operator with Dirichlet boundary condition. The reader can easily
check that the assumptions of Theorem \ref{esistenza} are
fulfilled and that the infimum is smaller than zero. Moreover the
infimum of the associated functional is achieved by a function of
the form $(u,0)$, or $(0,u)$, where $u\not\equiv0$; indeed, any
minimizer $(u_1,u_2)$ can be replaced by, say, $(u_1+u_2,0)$
keeping the same level of the functional. Hence the associated
variational problem does not have uniqueness of solutions.
\end{es}

A major obstruction to uniqueness is the lack of convexity that
may occur both in the Lagrangian (as the above example
illustrates) and in the constraint; nevertheless, the following
result shows that the full convexity of the Lagrangian is
sufficient to prove uniqueness of the minimizer, provided the
diffusions do not depend on $i$.
\begin{teo}\label{teo2} Under the assumptions of Section \ref{ipotesi},
assume moreover that
\neweq{uniq1} d_i\equiv d_j\;, \forall i,j\endeq
\neweq{uniq2}F_i(x,s) \mbox{ is concave in the variable $s$, for all $x\in\Omega$}
\endeq
 Then {\rm Problem \ref{Prb}} has an unique minimizer.
\end{teo}

\begin{proof} let $c=\inf\{J(U):\,U\in{\cal U}\}$.
Arguing by contradiction, we consider two minimizers $U=(u_1,\cdots,u_k)$
and $V=(v_1,\cdots,v_k)$ achieving $c$, with $u_i\not\equiv v_i$ for
some $i$. For every $i$ and $\l\in[0,1]$, we define
\begin{eqnarray}\label{hat1}
 \widehat u_i &=& u_i-\sum_{h\neq i}u_h\\\label{hat2}
 \widehat v_i &=& v_i-\sum_{h\neq i}v_h\\\nonumber
 w_i^{(\l)} &=& [\l\widehat u_i+(1-\l)\widehat v_i]^+.
 \end{eqnarray}
Our goal is to show that $J(w_i^{(\l)})< \l J(u_i)+(1-\l) J(v_i)$,
for every $\l\in(0,1)$. It is worthwhile noticing that this
property can be seen as a convexity type property combined with a
special type of projection on the constraint ${\cal U}$.\\ To
begin with, we have to show that the $ w_i^{(\l)}$'s satisfy the
constraint. We first notice that $w_i^{(\l)}\geq0$ and
$w_i|_{\partial\O}=\phi_i$. Furthermore we have
$$w_i^{(\l)}w_j^{(\l)}=0\mbox{ a.e. on }\Omega\;,\mbox{ for }j\neq
i.$$
Indeed, assume $w_i^{(\l)}(x)>0$; this means that
$$\l u_i(x)+(1-\l) v_i(x)>\sum_{h\neq i}\l u_h(x)+(1-\l) v_h(x)
\geq \l u_j(x)+(1-\l) v_j(x)\;\forall j\neq i.$$
Therefore we have, when $j\neq i$,

$$\l u_j(x)+(1-\l) v_j(x)<\l u_i(x)+(1-\l) v_i(x)\leq\sum_{h\neq j}\l u_h(x)+(1-\l) v_h(x)\;,$$

and hence $w_j^{(\l)}(x)\leq0$.

Let us denote the supports $$\Gamma_i^{(\l)}=\supp(w_i^{(\l)})\;;$$
recalling that $F_i(x,0)\equiv0$  we
note that
$$J\big(\sum_{i=1}^k
w_i^{(\lambda)}\big)=\sum_{i=1}^k\int_{\Gamma_i^{(\lambda)}}\big(\frac12d^2|\nabla
w_i^{(\l)}|^2 -F_i(w_i^{(\l)})\big)dx.$$
 In view of
\eq{uniq1}, using the convexity of the quadratic part of the
functional, and the definitions \eq{hat1},\eq{hat2} and keeping in
mind that both the $(u_i)'s$ and $(v_i)'s$ have disjoint supports,
we obtain that, for every $\l\in(0,1)$,
\begin{eqnarray*}
\sum_{i=1}^k&\int_{\Gamma_i^{(\lambda)}}{\big(\frac12 d^2|\nabla
w_i^{(\lambda)}|^2\big)dx}<\sum_{i=1}^k\int_{\Gamma_i^{(\lambda)}}\l\big(\frac12
d^2|\nabla \widehat u_i|^2\big)+(1-\l)\big(\frac12 d^2|\nabla
\widehat v_i|^2\big)dx\\&
\leq\sum_{i=1}^k\int_{\Omega}\l\big(\frac12d^2|\nabla
u_i|^2\big)+(1-\l)\big(\frac12d^2|\nabla  v_i|^2\big)dx
\end{eqnarray*}
Now we turn to the potential integral. By assumption \eq{uniq2}
and the evenness of the potentials $F_i$'s, the inequality
 $$-F_i(x,\l\widehat u_i(x)+(1-\l)\widehat
v_i(x))\leq-\l F_i(x,\widehat u_i(x))-(1-\l)F_i(x,\widehat
v_i(x))$$ holds whenever $x\in \supp (u_i)\cap\supp(v_i)$. Hence,
assume that $x\in \supp (u_i)\cap\supp(v_j)$, for some $j\neq i$;
let us fix $\l$ and let $x$ such that $\l u_i(x)-(1-\l)v_j(x)>0$
(the symmetric case is obtained by parity), and we   study the
term $-F_i(x,\l u_i(x)-(1-\l)v_j(x))$. Introducing the auxiliary
function $\Psi(\l)= -F_i(x,\l u_i-(1-\l)v_j)+ \l
F_i(x,u_i)+(1-\l)F_j(v_j)$, let $\bar\l\in(0,1)$ such that $\bar\l
u_i-(1-\bar\l)v_j=0$. It is easy to see that $\Psi(1)=0$, and
$\Psi''(\l)>0$. Moreover, $\Psi(\bar\l)\leq0$, since
$F_i(x,0)\equiv 0$ and,by \eq{uniq2}, $F_i\leq 0$. Hence, by
convexity, we infer that $\Psi(\l)<0$, for every $\l\in(\bar \l,1)$
and therefore $$-F_i(x,\left(\l \widehat u_i+(1-\l)\widehat
v_j\right)^+)\leq -\l F_i(x,u_i)-(1-\l)F_j(x,v_j)$$ holds
when $x\in\supp(u_i)\cap\supp(v_j)$.

Finally, gathering together all these inequalities, for every
fixed $\l$ we obtain $$J\big(\sum_{i=1}^k w_i^{(\lambda)}\big)<\l
J\big(\sum_{i=1}^k u_i\big)+(1-\l)J\big(\sum_{i=1}^k
v_i\big)<c\;,$$ a contradiction.
\end{proof}
The requirement of Theorem \ref{teo2}
that the Lagrangians are convex in all variables may seem
very restrictive.The following result makes a different assumption, still
sufficient for the uniqueness, that may be more useful in the applications,
for it is always satisfied (given the $d_i$'s and $F_i$'s) provided the
domain $\Omega$ is small enough.

\begin{coro} Let $d_i\in W^{2,\infty}$ be given.
Assume that $F_i$ are of class ${\cal C}^2$ in the variable $s$
and let us denote
\neweq{fsec}
b_i(x)= \sup_{s\in\R}\frac{\partial^2 F_i}{\partial s^2}(x,s)
\endeq
Assume that
there exists a positive function $d$ such that, for every $i=1,\dots,k$,
\neweq{supersol}
-\Delta d+\left(\frac{\Delta d_i}{d_i}-\frac{b_i}{2d_i^2}\right)d\geq 0\qquad\mbox{in $\Omega$}\;.
\endeq
Then,
{\rm Problem \ref{Prb}} has a unique minimizer.
\end{coro}
\begin{proof}
the following
identities hold:

\begin{eqnarray*}
\int_\Omega d^2_i(x)|\nabla u_i|^2
&=&\int_\Omega d^2(x)\left|\nabla\left(\frac{d_i u_i}{d}\right)\right|^2
-d^2\nabla\left(\frac{d_i}{d}\right)\cdot
\nabla\left(\frac{d_iu_i^2}{d}\right)\\&=&
\int_\Omega d^2(x)\left|\nabla\left(\frac{d_i u_i}{d}\right)\right|^2
+\div\left(d^2\nabla\left(\frac{d_i}{d}\right)\right)\frac{d_iu_i^2}{d}
-\int_{\partial\Omega}d d_iu_i^2\nabla\left(\frac{d_i}{d}\right)\cdot\nu
\end{eqnarray*}

Therefore, up to a constant (recall that our boundary data are
prescribed), by the change of variables $v_i=d_iu_i/d$ we can
transform the initial Lagrangians into ones of the following form:
$$
d^2(x)|\nabla v_i|^2+\frac{d}{d_i}
\div\left(d^2\nabla\left(\frac{d_i}{d}\right)\right)v_i^2-F_i(x,d
v_i/d_i).
$$
Recalling the definition of the $b_i$'s of \eq{fsec}, one easily checks that
these Lagrangians are convex in the variable $v_i$ provided
the following inequality holds for every index $i$
$$
2\frac{d}{d_i}
\div\left(d^2\nabla\left(\frac{d_i}{d}\right)\right)-\frac{d^2}{d^2_i}b_i\geq
0
$$
that is easily seen to be equivalent to \eq{supersol}.
\end{proof}

\section{Extremality conditions}\label{stoinesse}

The goal of this section is to prove that the minimizers of the variational
problem Problem \ref{Prb} satisfy a suitable set of differential inequalities.
To start with, let $(u_1,\dots,u_k)$ be a minimizer of problem \ref{Prb};
we define $\widehat
f(x,\widehat u_i)$
\neweq{hatf}
\widehat f(x,\widehat u_i)=\sum_{j} f_j(x,\widehat
u_i)\chi_{\supp(u_j)}=\left\{\begin{array}{ll}f_i(u_i)&\mbox{ if }
x\in\supp(u_i)\\-f_j(x,u_j)&\mbox{ if }x\in\supp(u_j),\; j\neq i.
\end{array}\right.
\endeq

Note that this definition is consistent with that of \eq{hat1},
for the functions $f_i$ are extended by oddness.
Our main result is the following:
\begin{teo}\label{teo1} Let
$U$ be a solution to {\rm Problem \ref{Prb}}. Then, for every $i$,
we have, in distributional sense,
\begin{enumerate}
\item[(i)] $-\Delta u_i\leq f_i(x,u_i)$ \item[(ii)] $-\Delta \widehat u_i\geq \widehat f(x,\widehat u_i)$ .
\end{enumerate}
\end{teo}

\begin{proof}
\begin{enumerate}
\item[(i)] We argue by contradiction. Then, there exists at least one index $j$ such that the claim
does not hold; that is, there is $0\leq\phi\in C^\infty_c(\Omega)$
such that
$$
\int_\Omega\nabla u_j\nabla\phi-f_j(x,u_j)\phi>0.
$$
For $0<t<1$ we define a new test function $V=(v_1,\dots,v_k)$ as
follows: $$ v_{i}=\left\{\begin{array}{ll}u_{i} &\mbox{
if $i\neq j$}
\\(u_i-t\phi)^+&\mbox{ if $i=j$.}\end{array}\right.$$
We claim that $V$ lowers the value of \eq{energia}; indeed we have
\begin{eqnarray*}J(V)-J(U)
&=&\int_{\Omega}\kern-3pt{\frac12}\left(|\nabla (u_j-t\phi)^+|^2-|\nabla u_j|^2\right)
-\int_{\Omega}\kern-3pt F_j(x,(u_j-t\phi)^+)-F_j(x,u_j) \\
&\leq&\kern-3pt\int_{\Omega}\kern-3pt\frac12\left(|\nabla u_j-t\phi|^2-|\nabla u_j|^2\right)
+t\int_{\Omega}\kern-3pt f_j(x,u_j)\phi +o(t) \\ &\leq& t\int_\Omega(-\nabla u_j\nabla\phi+f_j(x,u_j)\phi)+o(t).
\end{eqnarray*}
Choosing $t$ sufficiently small, we obtain $$J(V)-J(U)<0,$$ a contradiction.

\item[(ii)] Let $j$ and $0<\phi\in C^\infty_c(\Omega)$ such that
$$\int_\Omega\nabla \widehat u_j\nabla\phi -\widehat f(x,\widehat u_j)\phi\,dx<0.$$
Again, we show that the value of the functional can be lessen by
replacing $U$ with an appropriate new test function $V$. To this
aim we consider the positive and negative parts of $\widehat
u_j+t\phi$ and we notice that, obviously,
$$\{(\widehat u_j+t\phi)^->0\}\subset\{(\widehat u_j)^->0\}=\cup_{i\neq j}\supp(u_i)\;.$$
Let us define $V=(v_1,\dots,v_k)$ in the following way:
$$ v_i=\left\{
\begin{array}{ll}
\left(\widehat u_j+t\phi\right)^+, &\mbox{if $i=j$}\\
\left(\widehat u_j+t\phi\right)^-\chi_{\supp(u_i)}, &\mbox{if $i\neq j$}.
\end{array}\right.$$

We compute, using the definition \eq{hatf},

\begin{eqnarray*}
J(V)-J(U)&=&\sum_{i=1}^k\int_\Omega\frac12
\left(|\nabla v_i|^2-|\nabla u_i|^2\right)-
\left(F_i(x,v_i)-F_i(x,u_i)\right)\\
&=&\int_{\Omega}\frac12\left(|\nabla \widehat
u_j+t\phi|^2-|\nabla \widehat
u_j|^2\right)-\left(F_j(x,\left(\widehat u_j+t\phi\right)^+)-F_j(x,u_j)
\right)\\&-&\sum_{i\neq j}\left( F_i\left(x,\left(\widehat
u_j+t\phi\right)^-\chi_{\supp(u_i)}\right)-F_i(x,u_i)\right)
\\&=&
t\int_{\Omega}\nabla \widehat u_j\nabla\phi-
\int_{\Omega_j}
f_j(x, u_j)\chi_{\supp(u_j)}\phi+ \sum_{i\neq j}f_i(x,u_i)\chi_{\supp(u_i)}\phi+o(t)
\\ &=& t\int_{\Omega}\nabla \widehat u_j\nabla\phi-
\int_{\Omega}\widehat
f_j(x,\widehat u_j)\phi+o(t)\;.
\end{eqnarray*}
For $t$ small enough we find $J(V)<J(U)$, a contradiction.
\end{enumerate}
\end{proof}

\section{The class $\cal S$ and its basic properties}\label{sez.classe.esse}

Let $\cal U$ be the set of admissible $k$--uples as defined in Problem \ref{Prb} and let  $f_i$ be given satisfying $(A1),(A2)$;
we introduce the following functional class:

\begin{defin}\label{esse}
$$\cal S=\left\{
\begin{array}{l}
(u_1,\cdots,u_k)\in{\cal U}:
-\Delta u_i\leq f_i(x,u_i)\;,
 -\Delta \widehat u_i\geq \widehat f(x,\widehat u_i),\forall i=1,\dots,k
\end{array}\right\}
$$
\end{defin}

By virtue of Theorem \ref{teo1} the class $\cal S$ is the natural framework
where to develop our theory of regularity and free boundary of
minimizers to Problem \ref{Prb}. As already noticed in the introduction, the class $\cal S$
is of independent interest, for it contains the asymptotic limits of highly
competing diffusion systems.

Let us start with the following definitions:

\begin{defin}\label{def1}\rm
The multiplicity of a point $x\in\Omega$ is
$$m(x)=\sharp\left\{i:meas\left(\left\{u_i>0\right\}\cap B(x,r)
\right)>0 \;\forall \,r>0\right\}\;.$$ We shall denote by
$${\mathcal Z}_h(U)=\{x\in\Omega:m(x)\geq h\}$$
the set of points of multiplicity greater or equal than $h\in \bbbn$.
\end{defin}

The following properties are straightforward consequences of the
definition of $\cal S$ joint to the locally Lipschitz continuity of
the $f_i$'s that implies the validity of the Maximum Principle for
elliptic equations.

\begin{prop}\label{basic} Let $x\in\Omega$:
\begin{itemize}
\item[(a)] If $m(x)=0$, then there is $r>0$ such that $u_i\equiv0$ on $B(x,r)$, for every $i$.
\item[(b)] If $m(x)=1$, then there are $i$ and $r>0$ such that $u_i>0$  and
$-\Delta u_i=f_i(x,u_i)$ on $B(x,r)$.
\item[(c)] If $m(x)=2$, then are $i,j$ and $r>0$ such that $u_k\equiv 0$ for $k\neq i,j$ and
$-\Delta (u_i-u_j)=f_{i,j}(x,u_i-u_j)$ on $B(x,r)$, where $f_{i,j} (x,s)=f_i(x,s^+)-f_j(x,s^-)$.
\end{itemize}
\end{prop}
\proof part (a) follows directly by the definition of multiplicity;
assume $m(x)=1$ and let $r>0$ be such that $meas(\supp(u_j)\cap B(x,r))=0$ for every $j$ except
$i$. Then, by definition of $\cal S$, $u_i$ satisfies the equation $-\Delta u_i=f_i(x,u_i)$ on
$B(x,r)$. We can write
$$-\Delta u_i=a_i(x)u_i=\frac{f_i(x,u_i)}{u_i(x)}u_i$$
where $a_i\in L^\infty$ by the Lipschitz continuity of the $f_i$. Then, since $u_i\geq 0$, we infer from the
strong maximum principle that $u_i>0$ on $B(x,r)$. The second statement follows immediately
from the definition.
\endproof

\begin{remark}\label{deriv.norm}$ $\\
\rm  - We can not exclude, at this stage, the occurrence of points of multiplicity zero,
although this possibility will be ruled out at the end of Section \ref{free.boundary}, at least in
two dimensions, under a weak non degeneracy assumption. Note that
$\partial\{x\in\O:m(x)=0\}\subset {\cal Z}_3\cup \partial\O$.\\
- The second point of Proposition \ref{basic} says in particular that,
if $m(x)=2$
$$\lim_{y\to x\atop y\in \supp(u_i)}\nabla u_i(y)=-\lim_{y\to x\atop y\in \supp(u_j)}
\nabla u_j(y)\;.$$
If, by the way, the above limit is not zero, it follows that the set $\{x:m(x)=2\}$ is locally a
 ${\cal C}^1$ manifold of dimension $N-1$. Of course the above equality has to be
changed, in presence of variable diffusions, into
$$\lim_{y\to x\atop y\in \supp(u_i)}d_i(x)\nabla u_i(y)=-\lim_{y\to x\atop y\in \supp(u_j)}
d_j(x)\nabla u_j(y)\;.$$
\end{remark}

A major goal in the subsequent analysis will concern the
geometrical properties of the supports of the densities $u_i$ and
their common boundaries. As a first consequence of the Maximum Principle we can give a criterium for the
connectedness of the supports:
\begin{prop} Assume, for some $i$, that $u_i$ is continuous and
$\supp(u_i)\cap\partial\Omega$ is connected.
If
\neweq{firsteig}\frac{f_i(x,s)}{s}<\lambda_1(\Omega)\endeq
where $\lambda_1$ denotes the first eigenvalue of the Dirichlet operator with zero boundary
condition, then $\supp(u_i)$ is connected.
\end{prop}
\begin{proof}
assume not, then there is a connected component $A$ of $\supp(u_i)$ that
does not touch the boundary; on the other hand, $u_i>0$
satisfies the equation
$$-\Delta u_i=\left(\frac{f_i(u_i)}{u_i} \right)u_i\;$$
on $A$ with vanishing boundary trace. When testing the equation with $u_i$, by \eq{firsteig} we obtain a contradiction with the Poincar\'e inequality.
\end{proof}

It is worthwhile noticing that the condition \eq{firsteig} is always satisfied
on small domains. This can be useful in the local analysis of the solutions.
Another useful property of $\cal S$ is that its elements are uniformly bounded in the interior of $\Omega$ thanks to the next proposition.

\begin{prop}\label{prop.pinzamento}
 Let $(u_i)$ be an element of $\cal S$. Then the following hold
\begin{itemize}
\item[(i)]
 There are functions $(\Phi_i)\in W^{1,\infty}_{\rm loc}(\Omega)$ such that
\neweq{sub}
\begin{array}{rcll}
-\Delta \Phi_i &=& f_i(x,\Phi_i)&\mbox{in }\Omega\\
\Phi_i&=&u_i&\mbox{on }\partial\Omega\\
\Phi_i & \geq& u_i&\mbox{in }\Omega
\end{array}
\endeq

\item[(ii)]
There are functions $(\Psi_i)\in W^{1,\infty}_{\rm loc}(\Omega)$ such that
\neweq{super}
\begin{array}{rcll}
-\Delta \Psi_i &=& \widehat f(x,\widehat{u}_i)&\mbox{in }\Omega\\
\Psi_i&=&\widehat u_i&\mbox{on }\partial\Omega\\
\Psi_i & \leq& \widehat u_i&\mbox{in }\Omega
\end{array}
\endeq
(in particular, $\Psi_i^+\leq \widehat u_i^+=u_i$).
\end{itemize}
\end{prop}

\proof to prove the first assertion we apply the method of upper-lower
solutions: we need an ordered pair of functions $\alpha_i\leq\beta_i$ where $\alpha_i$ is  subsolution and $\beta_i$
 supersolution of problem \eq{sub}. We simply let $\alpha_i=u_i$  as lower solution;
 on the other hand we obtain a suitable $\beta_i$ by solving
$$
\begin{array}{rcll}
-\Delta \beta_i &=& b_i(x)\beta_i&\mbox{ in }\Omega\\
\beta_i&=&M+\phi_i&\mbox{ on }\partial\Omega
\end{array}
$$
for large constants $M$.
Notice that assumption (A2) implies the existence of arbitrarily
large positive functions $\beta_i$ satisfying the above problem. Furthermore,
since  $b(x)\beta_i\geq f_i(x,\beta_i)$, then
the $\beta_i$ are supersolutions to equation \eq{sub}.
Finally we get $\beta_i\geq \alpha_i$ by the maximum principle.\\
The proof of the second assertion  is trivial, since the boundedness of $\widehat{u}_i$ and the assumption (A1)
imply the existence of a solution for problem \eq{super}; the relation $\Psi_i \leq \widehat u_i$ then
follows by the maximum principle.\\
To conclude we observe that the regularity of the $\Phi_i$'s and $\Psi_i$'s follows by the standard regularity theory for elliptic equations and our assumptions on the
boundary data and the nonlinearities.
\endproof
\begin{remark}\label{equivS-S*}{\rm
As a consequence of the above proposition,
the components of each element $U\in{\cal S}$ are uniformly bounded on compact subsets of $\O$.
Then, recalling (A1),(A2), if $\omega\subset\subset\O$  there exists $M>0$ (depending only on $\omega$)
such that
\neweq{emme}
U\in{\cal S}\;\;\;\Longrightarrow\;\;\;-\Delta u_i\leq M,\;\;\;-\Delta \widehat u_i\geq-M\hspace{1cm}\mbox{ on }\omega
\endeq
Furthermore, the regularity can be improved up to the boundary of
$\Omega$ in the sense that to bounded boundary data there
correspond bounded barriers $\Phi_i$ and $\Psi_i$. Moreover the
barriers will be Lipschitz continuous up to the boundary when both
the data and the boundary $\partial\Omega$ enjoy the same
regularity.

 }\

\end{remark}

 \section{The class ${\cal S}^*$}\label{sez.classe.esse*}

Let $M\geq 0$ and $h$ be a fixed integer. We introduce
\begin{defin}\label{esse*}
$${\cal S}^*_{M,h}(\omega)=\left\{(u_1,\cdots,u_h)\in
(H^1(\omega))^h:\,\begin{array}{l}
 u_i\geq0,\, u_i\cdot u_j=0\mbox{ if }i\neq j\\
 -\Delta u_i\leq M,\,-\Delta \widehat u_i\geq-M
\end{array}\right\}
$$
\end{defin}
It follows from Remark \ref{equivS-S*}, that ${\cal S} \subset {\cal S}^*_{M,h}$.
It will be more convenient to work in this larger class rather that
in the  class $\cal S$,
for it is closed with respect to the limits of sequences of scaled functions.
This property will be extremely useful in performing the blow--up analysis in Section \ref{regularity}.
We first present a technical result concerning with the elements
of ${\cal S}^*_{M,h}(B(0,1))$ when 0 is a point of multiplicity at
least two.
\begin{prop}\label{poincare.convergenze}
There exists $M^*$ such that, for all $0<M<M^*$ the following
holds
\begin{itemize}
\item[(a)] Let $U_n\in  {\cal S}^*_{M,h}(B(0,1))$ such that
$m(U_n)(0)\geq 2$ and $U_n\rightharpoonup U$ with $U\not\equiv 0$.
Then $$\sharp\{i=1,\dots,h\;:\:u_i\not\equiv 0\}\geq 2.$$
\item[(b)] For every $\gamma>0$ there exists $C(\gamma)>0$ such that, if $U\in  {\cal S}^*_{M,h}(B(0,1))$ with
$m(U)(0)\geq 2$ and $\|U\|_{L^2(B(0,1))}\geq \gamma$, then  $$\int_{B(0,1)}|\nabla
U|^2\geq C(\gamma) \int_{B(0,1)}|U|^2$$ where $C>0$ only depends on $M^*$ and $\gamma$.
\item[(c)] Let $U_n\in  {\cal S}^*_{M,h}(B(0,1))$ such that
$m(U_n)(0)\geq 2$; if $\|\nabla U_n\|_{L^2(B(0,1))}=1$ and
$\|\nabla U_n\|_{L^2(\partial B(0,1))}$ is bounded, then there
exists $U\not\equiv 0$ such that $U_n\rightharpoonup U$.
\end{itemize}
\end{prop}
\proof {\it Part (a).} Arguing by contradiction we can assume the
existence of $M_n\to 0$ and $U_n\in {\cal S}^*_{M_n,h}(B(0,1))$
such that $m(U_n)(0)\geq 2$, $U_n\rightharpoonup U$ with
$U\not\equiv 0$ and  $N:=\sharp\{i=1,\dots,h\;:\:u_i\not\equiv
0\}\leq 1$. Note that, since $U\not\equiv 0$, then $N$ must be
equal to 1. Thus  we can assume that $U$ is of the form
$U=(u_1,0,\dots,0)$ where $u_1\not\equiv 0$.
 Now
consider $\phi>0$ solution of the problem $-\Delta\phi=1$ on
$B(0,1)$ with boundary conditions $\phi=0$. It holds
$$-\Delta(\widehat{u}_{1,n}+M_n\phi)\geq 0$$ and
by the mean value property for superharmonic functions we have
\neweq{mean.value}
0=\widehat{u}_{1,n}(0)\geq
\int_{B(0,1)}\widehat{u}_{1,n}(x)dx+M_n\left(\int_{B(0,1)}\phi(x)dx-\phi(0)\right).
\endeq
By the compact embedding of $H^1$ in $L^1$ and since
$\widehat{u}_1=u_1$, it holds
$$\int_{B(0,1)}\widehat{u}_{1,n}(x)dx\to \int_{B(0,1)}u_1(x)dx=\a>0.$$
On the other hand, the last term of \eq{mean.value} becomes less
than $\a/2$ when $M_n\to 0$: this finally provides the
contradiction $0>\a/2$.

{\it Part (b).} Let $\gamma>0$ be fixed. Assume by contradiction the existence of $M_n\to
0$ and $U_n\in {\cal S}^*_{M_n,h}(B(0,1))$ such that
$m(U_n)(0)\geq 2$, $\|U_n\|_{L^2(B(0,1))}\geq \gamma$ and
$$\frac{\int_{B(0,1)}|\nabla U_n|^2}{\int_{B(0,1)}|U_n|^2}\to 0.$$
Let us define $V_n:=U_n\cdot\|U_n\|^{-1}_{L^2(B(0,1))}$: note that
$V_n\in {\cal S}^*_{M_n/\gamma,h}(B(0,1))$, $m(V_n)(0)\geq 2$; furthermore
$\|V_n\|_{L^2(B(0,1))}= 1$ and $\|\nabla V_n\|_{L^2(B(0,1))}\to 0$
by construction. Then there exists $V$ such that
$V_n\rightharpoonup V$; moreover $V\not\equiv 0$, since by the
compact embedding of $H^1$ in $L^2$ it holds $\|V\|_{L^2(B(0,1))}=
1$. Now note that, since the gradients vanish in the $L^2$ norm,
then $V=(c_1,\dots,c_h)$ where $c_i\in \bbbr$. Furthermore, since
the supports of the components of $V_n$ are mutually disjoint,
passing to the limit a.e. in $B(0,1)$, we get $c_i\cdot c_j=0$ if
$i\neq j$: this means that only one of the components of $V$ is
not identically 0. This is in contradiction with Part (a) when
applied to the
sequence $V_n$.\\

{\it Part (c).} By  assumption $\|\nabla U_n\|_{L^2(B(0,1))}$ is
bounded: then, if $M$ is small enough we can apply Part (b) and
thus the whole $H^1$--norm  $\| U\|_{H^1(B(0,1))}$ is bounded. This
provides $U$ such that $U_n\rightharpoonup U$. Now test the
variational inequality $-\Delta u_{i,n}\leq M$ with $u_{i,n}$:
$$\int_{B(0,1)}|\nabla u_{i,n}|^2\leq
M \int_{B(0,1)} u_{i,n}+\int_{\partial B(0,1)}u_{i,n}
\frac{\partial}{\partial \nu}u_{i,n}.$$ Assume by
contradiction that $U\equiv 0$: by the compact embedding of $H^1$
in $L^p$ and since $\| \nabla U_n\|_{L^2(\partial B(0,1))}$ is
bounded, we deduce that the r.h.s vanishes. This implies $\|\nabla
U_n\|_{L^2(B(0,1))}\to 0$, in contradiction with the assumption
$\|\nabla U_n\|_{L^2(B(0,1))}=1$.
\endproof
\begin{rem}\label{proprieta_media}{\rm The argument used in the proof of part (a) allows
to establish a mean value property for functions which are
superharmonic up to a small term. Precisely, if $(v_n)\subset H^1$
is such that $-\Delta v_n\geq -M_n$ with $(0<)M_n\to 0$ and
$\int_{B(0,1)}v_n\geq \a >0$, then $v_n(0)\geq \a/2$ if $n$ is large enough. }
\end{rem}
\section{Lipschitz Regularity}\label{regularity}

A key tool in studying the regularity of both the function $U$ and the free boundary is
a suitable version of the celebrated {\it monotonicity theorem}, see \cite{acf,cks}.
In this paper we shall take advantage of the following formula which is proven in
\cite{cjk}.
\begin{lemma}\label{monot.limit}
Let $w_1,w_2\in H^1\cap L^\infty$ such that $-\Delta w_i\leq 1$, $w_1\cdot w_2=0$ a.e. and
$x_0\in\partial(\supp(w_i))$, $i=1,2$. Then there exists $C>0$,
independent of $x_0$, such that
\neweq{monot2}
\prod_{i=1}^{2}{1\over r^2}\int_{B(x_0,r)}\frac{|\nabla
w_i(x)|^2}{|x-x_0|^{N-2}}dx \leq C.
\endeq
\end{lemma}

In the proof of our regularity results we shall also need the following technical lemma.
\begin{lemma}\label{cambia(x,r)}
Let $U\in H^1(\O)$ and let us define
\begin{equation}\label{def.phi} \phi(x,r):=\frac1{r^N}\int_{B(x,r)
\cap\O}|\nabla U(y)|^2dy.
\end{equation}
If $(x_n,r_n)$ is a sequence in
$\overline\O\times\bbbr^+$ such that
$$
\phi(x_n,r_n)\to\infty
$$
 then $r_n\to0$ and
\begin{itemize}
\item[(i)] there exists a sequence $(r'_n)_n\subset\bbbr^+$
such that $\phi(x_n,r'_n)\to\infty$,  and
$$\int_{\partial
B(x_n,r'_n)\cap\O}|\nabla U|^2\leq {N\over r'_n}
\int_{B(x_n,r'_n)\cap\O}|\nabla U|^2;$$
\item[(ii)] if $A\subset\overline\O$ and
$$
\frac{dist(x_n,A)}{r_n}\leq C
$$
then there exists a sequence $(x'_n,r'_n)$ such that
$\phi(x'_n,r'_n)\to\infty$ and $x'_n\in A$ for every $n$.
\end{itemize}
\end{lemma}

\proof since $U\in H^1(\O)$ and $\phi(x_n,r_n)\to+\infty$,
obviously $r_n\to0$. We begin proving (i). Let $g$ be defined on
the whole $\bbbr^N$ as
$$
g(x):=\left\{\begin{array}{cl}
      |\nabla U(x)|^2 & x\in\O\\
      0 & x\in\bbbr^N\setminus\O.
      \end{array}\right.
$$
Clearly $\int_{B(x,r)\cap\O}|\nabla U|^2=\int_{B(x,r)}g$,
$\int_{\partial B (x,r)\cap\O}|\nabla U|^2=\int_{\partial
B(x,r)}g$. We observe that
$$
\frac{\partial \phi}{\partial r}(x_n,r_n)=\frac{1}{r_n^N}
\left(\int_{\partial B(x_n,r_n)}g(y)dy-\frac{N}{r_n}
\int_{B(x_n,r_n)}g(y)dy\right).
$$
As a consequence, our problem is reduced to find $r'_n$ such that
$\phi_r(x_n,r'_n)\leq0$ and $\phi(x_n,r'_n)\to+\infty$. Let
$r'_n:=\inf\{r\geq r_n:\,\phi_r(x_n,r)\leq0\}$. We have that
$r'_n<+\infty$ for every $n$ (recall that
$\phi(x_n,r)\to0<\phi(x_n,r_n)$ as $r\to\infty$), and
$\phi(x_n,r_n)\leq\phi(x_n,r'_n)$, and (i) is proved. In order to
prove (ii), let $x'_n\in A$ such that $dist(x_n,x'_n)\leq
2Cr_n$ (this is possible by assumption) and let $r'_n:=(2C+1)r_n$. By construction
$B(x'_n,r'_n)\supset B(x_n,r_n)$, that implies
$\phi(x'_n,r'_n)\geq(2C+1)^{-N}\phi(x_n,r_n)$, and also (ii)
follows.
\endproof

\subsection{Local Lipschitz continuity in ${\cal S}^*$}
The class ${\cal S}^*$ seems to be the natural framework for proving
the interior Lipschitz regularity. Indeed we have:
\begin{teo}\label{teo.lip.interno}
Let $M>0$ and $k$ be a fixed integer. Let $U\in\cal S^*_ {M,k}(\O)$: then $U$ is
Lipschitz continuous in the interior of $\Omega$.
\end{teo}

\begin{proof}
we consider the set $\Omega'\subset\subset\Omega$ compactly
enclosed in $\Omega$ and the function $\phi(x,r)$ defined as in
\eq{def.phi}, restricted to the set $D:=\{x\in\O',
r\in\bbbr^+:2r<dist(\partial\O',\partial\O)\}$. We have to prove
that $\phi$ is bounded on $D$.
 We argue by contradiction, assuming that
$$\sup_D\phi(x,r)=+\infty$$ i.e. there exists a
sequence $(x_n,r_n)$,  such that
\begin{equation}\label{eq1}
\lim_{n\to\infty}\frac1{r_n^N}\int_{B(x_n,r_n)}|\nabla
U(y)|^2dy=+\infty.\end{equation} By Lemma \ref{cambia(x,r)},(i),
there exists a sequence (denoted again by $r_n$) satisfying
\eq{eq1} and moreover
\neweq{fifi}\int_{\partial B(x_n,r_n)}|\nabla U|^2\leq {N\over r_n}
\int_{B(x_n,r_n)}|\nabla U|^2.
\endeq
Now we distinguish three cases, according to the nature of $x_n$
and $r_n$, up to suitable subsequences:

CASE I: $m(x_n)=0$ for all $n$ (up to a subsequence) and
$\frac{dist(x_n,{\cal Z}_1)}{r_n}\geq 1$.\\
We immediately obtain a contradiction with \eq{eq1}, since in this
case $u_i\equiv 0$ on $B(x_n,r_n)$, for all $i$. As a consequence,
using Lemma \ref{cambia(x,r)},(ii),  without loss of
generality, we can assume that $m(x_n)\geq1$ for every $n$ (besides \eq{eq1} and \eq{fifi}).

CASE II: $m(x_n)=1$ for all $n$ and
$\frac{dist(x_n,{\cal Z}_2)}{r_n}\geq 1$.\\
In this case we can assume that only $u_1\not\equiv 0$ on $B(x_n,r_n)$ and we
define
$$v_n(x)=(u_1(x)-u_1(x_n))^+\hspace{1cm}w_n(x)=(u_1(x)-u_1(x_n))^-$$
Then $v_n$ and $w_n$ are in $L^\infty(\O')$ by Remark \ref{esse}, they have disjoint  supports, $-\Delta v_n\leq M$, $-\Delta w_n\leq M$ and
$v_n(x_n)=w_n(x)=0$.
Thus we can apply the
monotonicity formula of Lemma \ref{monot.limit}
$$
\frac{1}{r_n^2}\int_{B(x_n,r_n)}\frac{|\nabla
v_n(x)|^2}{|x-x_n|^{N-2}}\cdot\frac{1}{r_n^2}\int_{B(x_n,r_n)}\frac{|\nabla
w_n(x)|^2}{|x-x_n|^{N-2}}\leq C
$$
 uniformly. Since
$|x-x_n|<r_n$, we deduce $$\frac1{r_n^N}\int_{B(x_n,r_n)}|\nabla
v_n(x)|^2\cdot \frac1{r_n^N}\int_{B(x_n,r_n)}|\nabla
w_n(x)|^2\leq C $$where $C$ is independent of $n$. Comparing with \eq{eq1} we
have that only one of the two term is unbounded and forces the second one to vanish, e.g.
\neweq{awto0}\frac1{r_n^N}\int_{B(x_n,r_n)}|\nabla
v_n|^2\longrightarrow\infty\hspace{1cm}\frac1{r_n^N}\int_{B(x_n,r_n)}|\nabla
w_n|^2\longrightarrow0\endeq as $n\to\infty$.
Let $$L_n^2:=\frac1{r_n^N}\int_{B(x_n,r_n)}|\nabla
v_n|^2$$ and let us perform the blow up analysis around $x_n$ with parameter $L_n$
 by defining $V_n=(v_{1,n},v_{2,n})$ as
$$v_{1,n}(x)=\frac1{L_nr_n}v_n(x_n+r_nx)\hspace{1cm}
v_{2,n}(x)=\frac1{L_nr_n}w_n(x_n+r_nx)\hspace{1cm}x\in B(0,1).$$
It is easy to verify that $V_n\in {\cal S}^*_{M_n,2}(B(0,1))$,
where $M_n=r_nM/L_n\to 0$. By construction we have that
$\int_{B(0,1)}|\nabla V_n|^2$ is bounded; using \eq{fifi}, this implies that $\int_{\partial
B(0,1)}|\nabla V_n|^2$ is bounded too: thus $V_n$ satisfies the
assumptions of Proposition \ref{poincare.convergenze}. This
provides the existence of a weak limit $V=(v_1,v_2)$ (by Part (c))
such that $v_i\not\equiv 0$, $i=1,2$ (by Part (a)). But this is in
contradiction with \eq{awto0} that forces $v_2\equiv 0$ (by Part
(b)). Again, this contradiction and Lemma \ref{cambia(x,r)},(ii)
allow us to assume $m(x_n)\geq 2$.

CASE III: $m(x_n)\geq 2$ for all $n$.\\
In this case the proof follows closely the line of the previous
one. Let us give some details: since $m(x_n)\geq 2$, we can apply the
monotonicity formula to each pair
$w_1:=u_l$, $w_2:= u_k$:
$$
\frac{1}{r_n^2}\int_{B(x_n,r_n)}\frac{|\nabla
w_1(x)|^2}{|x-x_n|^{N-2}}\cdot\frac{1}{r_n^2}\int_{B(x_n,r_n)}\frac{|\nabla
w_2(x)|^2}{|x-x_n|^{N-2}}\leq C
$$
 uniformly. By \eq{eq1} we deduce the existence of one index $i$ such that, up to a
 subsequence, it holds
\neweq{wto0}\frac1{r_n^N}\int_{B(x_n,r_n)}|\nabla
u_i|^2\longrightarrow\infty\hspace{1cm}\frac1{r_n^N}\int_{B(x_n,r_n)}|\nabla
u_j|^2\longrightarrow0\hspace{1cm} \forall j\neq i\endeq as $n\to\infty$.
Let $$L_n^2:=\frac1{r_n^N}\int_{B(x_n,r_n)}|\nabla
u_i|^2$$ and perform the blow up analysis around $x_n$ with parameter $L_n$
 by defining $V_n=(v_{j,n})$ as
$$v_{j,n}(x)=\frac1{L_nr_n}u_j(x_n+r_nx)\hspace{1cm}x\in B(0,1).$$
It is easy to verify that $V_n\in {\cal S}^*_{M_n,k}(B(0,1))$,
where, again, $M_n=r_nM/L_n\to 0$. By construction we have that
$\int_{B(0,1)}|\nabla V_n|^2$ is bounded and, again by \eq{fifi}, $\int_{\partial B(0,1)}|\nabla V_n|^2$ is
bounded too: thus $V_n$ satisfies the assumptions of Proposition
\ref{poincare.convergenze}. This provides the existence of a weak
limit $V$ (by Part (c)) such that at least two of its components
are strictly positive (by Part (a)). But this is in contradiction
with \eq{wto0} that forces (by Part (b)) $v_j\equiv 0$ for all
$j\neq i$.
\end{proof}

\subsection{Lipschitz continuity up to the boundary in ${\cal S}$}

In this Section we are concerned with the regularity of the
elements of ${\cal S}$ up to the boundary, in the case of regular boundary and
Lipschitz boundary data. Our main goal is the following result.

\begin{teo}\label{teo.lip.globale}
Let $\partial\Omega$ be of class $C^1$, $U\in{\cal S}$ with
$u_i|_{\partial\O}=\phi_i$ and $\phi_i\in
W^{1,\infty}(\partial\O)$ for every $i$. Then $U\in
W^{1,\infty}(\overline\Omega)$.
\end{teo}
The proof relies upon the local analysis as developed in the previous section, joint
with the suitable use of the {\it pinching} property stated in Proposition
\ref{prop.pinzamento}. We begin with some preliminary remarks.

\begin{rem}\label{remlipglob}
Under the assumptions of Theorem \ref{teo.lip.globale} we have
\begin{itemize}
\item[(i)] $u_i\in C(\overline\O)$ for every $i$ (and, in particular,
it makes sense to consider pointwise values of $u_i$);
\item[(ii)] $\frac{\partial{u_i}}{\partial \nu}|_{\partial\O}\in
L^\infty(\partial\O)$ for every $i$.
\end{itemize}
\end{rem}

\proof since $U\in{\cal S}$, through Proposition
\ref{prop.pinzamento} we obtain the existence of $k$--uples of
functions $(\Phi_i)$, $(\Psi_i)$, with the properties introduced
in that proposition. Moreover, since $\phi_i\in
W^{1,\infty}(\partial\O)$ for every $i$, by standard regularity
theory for elliptic equations we infer $\Phi_i\in
W^{1,\infty}(\overline\Omega)$, $\Psi_i\in
W^{1,\infty}(\overline\Omega)$ for every $i$. By Theorem
\ref{teo.lip.interno} (which holds in this case through Remark
\ref{equivS-S*}) $u_i\in C(\O)$; since $\Psi_i^+\leq u_i
\leq\Phi_i$, and $\Psi_i^+\equiv\Phi_i$ on $\partial\O$, (i)
easily follows. Moreover, using the same inequality and the very
definition of directional derivative, we obtain, in distributional
sense,
$$
\frac{\partial \Psi_i^+}{\partial \nu}\geq \frac{\partial
u_i}{\partial \nu}\geq \frac{\partial \Phi_i}{\partial \nu},
$$
and also (ii) follows.\endproof

Now we are ready to prove Theorem \ref{teo.lip.globale}.

{\bf Proof of Theorem \ref{teo.lip.globale}:} let $U$ satisfy the
assumptions of the theorem. By the first part of the previous
remark $u_i\in L^\infty(\O)$ for every $i$, and hence there exists
a constant $M$ such that
$$
-\Delta u_i\leq f(x,u_i)\leq M.
$$
Moreover, by Proposition \ref{prop.pinzamento}, there exist
$k$--uples $(\Psi_i)$, $(\Phi_i)$ such that $\Psi_i^+(x)\leq u_i
(x)\leq\Phi_i(x)$ on $\overline\O$, $\Psi_i^+(x)=u_i(x)=\Phi_i(x)$
on $\partial\O$ and, as we just observed, $\Psi_i$, $\Phi_i\in
W^{1,\infty}(\overline\O)$ for every $i$.

As usually we define $\phi(x,r)$ as in \eq{def.phi} and we assume
by contradiction the existence of a sequence $(x_n,r_n)$ such
that $\phi(x_n,r_n)\to+\infty$. Assume that (up to a subsequence)
$\frac{dist(x_n,\partial\O)}{r_n}\geq 1$. This means that
$B(x_n,r_n)\cap\partial\Omega=\emptyset$ for every $n$. In this
situation, one can repeat exactly the same proof of Theorem
\ref{teo.lip.interno} (roughly speaking, in such a situation the
blow--up procedure does not ``see'' the boundary), obtaining the
same contradictions. Hence we can assume that $\frac{dist(x_n,
\partial\O)}{r_n}\leq 1$ and, by Lemma \ref{cambia(x,r)},(ii), without
loss of generality we can choose $x_n\in\partial\O$; moreover, we
take $r_n$ such that the inequality in Lemma
\ref{cambia(x,r)},(i) holds. Let $ \phi(x_n,r_n)=:L_n^2\to+\infty$
and define, for every $i$,
$$
u_{i,n}(x):=\frac{1}{r_n L_n}(u_i(x_n+r_nx)-u_i(x_n)),
$$
$$
\Psi_{i,n}(x):=\frac{1}{r_n L_n}(\Psi^+_i(x_n+r_nx)-
\Psi_i^+(x_n)),
$$
$$
\Phi_{i,n}(x):=\frac{1}{r_n L_n}(\Phi_i(x_n+r_nx)-\Phi_i(x_n)).
$$
We have that $\Psi_{i,n}(x)\leq u_{i,n}(x)\leq\Phi_{i,n}(x)$ on
$\overline\O_n$, where $\O_n:=\{x\in B(0,1):\,x_n+r_nx\in\O\}$
(this inequality holds because the non-scaled functions coincide
in $x_n\in\partial\O$). We observe that, taking into account that
$r_n\to0$ and $L_n\to\infty$, we have
$$
-\Delta u_{i,n}(x)=-\frac{r_n}{L_n}\Delta u_i(x_n+r_nx)\leq
\frac{r_n}{L_n}M\leq1
$$
for every $i$, when $n$ is sufficiently large. Testing the above inequality on $u_{i,n}$ we obtain
$$
\int_{\O_n}|\nabla u_{i,n}|^2\leq\int_{\O_n}u_{i,n}+\int_{\partial
\O_n}u_{i,n}\cdot\frac{\partial u_{i,n}}{\partial\nu}\leq\|u_{i,n}
\|_{L^\infty(\O_n)}+\|u_{i,n} \|_{L^2(\partial\O_n)}\left\|
\frac{\partial u_{i,n}}{\partial\nu}\right\|_{L^2(\partial\O_n)}.
$$
Our aim is to prove that the righthand side of the previous
inequality tends to 0 for every $i$. This will provide a
contradiction with the fact that, by construction, $\sum\int_{
\O_n}|\nabla u_{i,n}|^2=1$ for every $n$.

Clearly $\Psi_{i,n}(0)=\Phi_{i,n}(0)=0$. Moreover,
$$
\|\nabla \Phi_{i,n}\|_{L^\infty(\O_n)}=\frac{1}{L_n}\|\nabla
\Phi_{i}\|_{L^\infty(\O\cap B(x_n,r_n))}\to0,
$$
and the same holds for $\Psi_{i,n}$. This implies
$\|\Phi_{i,n}\|_{W^{1,\infty}(\overline\O_n)}\to0$,
$\|\Psi_{i,n}\|_{W^{1,\infty}(\overline\O_n)}\to0$ and therefore
$$
\|u_{i,n}\|_{L^\infty(\O_n)}\to0\quad\forall i.
$$
Since $\partial\O$ is of class $C^1$, we obtain that
$\partial\O_n$ has bounded $(N-1)$--dimensional measure, and thus
also
$$
\|u_{i,n}\|_{L^2(\partial\O_n)}\to0\quad\forall i.
$$
Therefore the only thing that remains to prove is that
$\|\frac{\partial u_{i,n}}{\partial
\nu}\|_{L^2(\partial\O_n)}$ is bounded. To this aim, let
$\partial\O_n=\Gamma_{1,n}\cup\Gamma_{2,n}$, where $\Gamma_{1,n}:=
\partial B(0,1)\cap\overline\O_n$ and $\Gamma_{2,n}:=\partial\O_n
\setminus\Gamma_{1,n}$. Then the estimate on $\Gamma_{1,n}$
descends from Lemma \ref{cambia(x,r)},(i), recalling that it implies
$$
\int_{\Gamma_{1,n}}\left|\frac{\partial u_{i,n}}{\partial
\nu}\right|^2\leq \frac{1}{r_n^{N-1}L_n^2}\int_{\partial B(x_n,r_n)\cap\O}|\nabla u_i|^2
\leq \frac{N}{r_n^{N}L_n^2}\int_{B(x_n,r_n)\cap\O}|\nabla u_i|^2=N.
$$
On the other hand, the estimate on $\Gamma_{2,n}$ is an easy
consequence of the bounded measure of $\Gamma_{2,n}$ and of Remark
\ref{remlipglob},(ii).
\endproof

\section{Further regularity in dimension N=2}\label{N=2}
\subsection{Vanishing of the gradient at multiple points}\label{gradientenullo}

Let $N=2$ and $U\in\cal S$. The main goal of this section is to
prove that the gradient of $U$ vanishes continuously at points of
multiplicity at least three. This result will be established
through the application of a monotonicity formula with three or
more phases. To start with, we need the following technical
result, that allows us to reduce the $u_i$'s to solutions to
suitable divergence--type equations.
\begin{lemma}\label{ker.div}
Let $a\in L^\infty$ and let $v$ be an $H^1$ solution of $-\Delta
v\leq a(x) v$ in $B(x_0,r)$. Then, if $r$ is small enough, there
exists $\f\in C^1$ such that $\f$ is radial with respect to $x_0$,
$\inf_{B(x_0,r)} \f>0$ and
$$-\div(\f^2(\nabla\frac{v}{\f}))\leq 0\hspace {1cm}\mbox{ in
}B(x_0,r).$$
\end{lemma}
\proof let us consider the eigenvalue problem
$$\left\{\begin{array}{lll} -\Delta u(x)=a(x)u(x)& x\in B(x_0,r)\\
u(x)>0 & x\in B(x_0,r)
\end{array}
\right. $$ If $r$ is small enough, the above problem can be solved
in the class of radial functions with respect to $x_0$. Let $\f$
be such a solution: now  by elementary computations $$-\f^2
\Delta\frac{v}{\f}-2\f\nabla\frac{v}{\f}\nabla\f\leq 0$$ giving
the required inequality. \endproof
This local reduction will be widely exploited throughout the present and the next section.
As a first application it allows to prove a variant of the original
monotonicity formula by Alt--Caffarelli--Friedman \cite{acf}.

\begin{lemma}\label{monot.limit.n=2}
Let $U\in \cal S $ and $w_i=\sum_{j\in I_i}u_j$, where
$I_1\cup...\cup I_h\subset\{1,...,k\}$. Assume that
$x_0\in\partial(\supp(w_i))$. Then for all $h\geq 2$ there exists $C$, independent of $x_0$
such that
\neweq{monot3}
\prod_{i=1}^{h}{1\over r^h}\int_{B(x_0,r)}|\nabla w_i(x)|^2dx \leq
C.
\endeq
\end{lemma}
\proof(sketch). Let $x_0\in\partial\supp(u_i)$ and recall that,
since $U\in\cal S$, it holds $-\Delta u_i\leq f_i(u_i)$ for all
$i$. Let $a(x):=\max_{i=1,\dots k}
\left|\frac{f_i(u_i(x))}{u_i(x)}\right|$ and note that $a\in
L^\infty$ by the assumption on $f_i$. By Lemma \ref{ker.div},
there exists $r>0$ and a regular radial function $\f$ which is
strictly positive on $B(x_0,r)$ such that, for all $i$
$$-\div(\f^2(\nabla\frac{u_i}{\f}))\leq 0\hspace {1cm}\mbox{ in
}B(x_0,r).$$ Now consider $w_i=\sum_{j\in I_i}u_j$, where
$I_1\cup...\cup I_h\subset \{1,...,k\}$ and let $\tilde
w_i:=w_i/\f$: then $-\div(\f^2(\nabla\tilde{w}_i))\leq 0$. Let
$N=2$ and set
$$\Phi(r)= \prod_{i=1}^{h}{1\over
r^h}\int_{B(x_0,r)}\f^2(x)|\nabla \tilde{w}_i(x)|^2 dx.$$
Following the proof of  Lemma 5.1 in \cite{acf} we can compute
$$\Phi'(r)\geq h\frac{\Phi(r)}{r}\left(-h+\frac{2}{h}\inf_{\cal P}\sum_{i=1}^h \sqrt{\Lambda_i}\right)$$
where ${\cal P}=\{(\G_i)_{i=1}^h, \cup \G_i=\partial B(0,1), meas(\G_i\cap\G_j)=0\}$ and
$$\Lambda_i=\inf_{v\in H^1_0(\Gamma_i)} \frac{\int_{\G_i}|\partial_\theta
v|^2}{\int_{\G_i}|v|^2}.$$
Since $\Lambda_i\geq\left(\frac{\pi}{|\Gamma_i|}\right)^2$
we immediately obtain that $\Phi'(r)\geq 0$ in
$[0,r']$.\\ By this formula and since there exist positive $a<b$
such that $a<\f(x)<b$ for all $x\in B(x_0,r')$, it follows in
particular that \eq{monot3} holds.
\endproof

Now we are ready to prove the main result of this section.

\begin{teo}\label{grad.nullo}
 If  $x_0\in{\cal Z}_3$ then $|\nabla
U(x)|\to 0$ as $x\to x_0$.
\end{teo}

\proof assume by contradiction the existence of $x_n\to x_0$ and
$r_n\to 0$ such that
\neweq{ass}\lim_{n\to+\infty}{1\over r_n^2}\int_{B(x_n,r_n)}|\nabla
U|^2=\alpha\endeq $\alpha\in(0,\infty]$.

We first claim that equation \eq{ass} holds for $x_n\in {\cal Z}_3$.
The proof of this fact can be done as follows: let $\rho,r>0$ be fixed and let
$A_{\rho,r}=\{x\in B(\rho,x_0)\;:\;d(x,{\cal Z}_3(U))\geq r\}$. By
Proposition \ref{basic}, in $A_{\rho,r}$ we can give alternate
positive and negative sign to the $u_i$'s in such a way that the
resulting function $v$ locally solves an equation of the form
$-\Delta v=f(x,v)$ where $f(x,v(x))=f_i(x,v(x))$ if $v(x)=\pm
u_i(x)$. For $x\in A_{\rho,r}$ we define
$$\Phi(x)={1\over
r^2}\int_{B(x,r)}|\nabla U(y)|^2dy \;={1\over
r^2}\int_{B(0,r)}|\nabla U(x+y)|^2dy.$$ By elementary computations
 $-\Delta (|\nabla v(x)|^2)=-2|\nabla
v(x)|^2-2f_x(x,v)|\nabla v(x)|-2f_s(x,v)))|\nabla v(x)|^2$:
recalling that $v$ is bounded by Remark \ref{equivS-S*} and since $|\nabla U|=|\nabla v|$, this gives
$-\Delta \Phi\leq a\Phi$ on $A_{\rho,r}$, for some positive
constant $a$ (depending on $r$). Then, by an extension of the
maximum principle, there exists $C>0$ (independent of $r$) such
that $\max_{A_{r,\rho}}\Phi\leq C\max_{\partial A_r}\Phi$.
 This implies that
 for $n$ large enough, \eq{ass} holds with $\alpha=\alpha/2$ and
for a choice of $x'_n$ such that $d(x'_n,{\cal Z}_3(U))\asymp
r_n$; then taking $z_n\in{\cal Z}_3(U)$ such that
$|z_n-x'_n|=d(x'_n, {\cal Z}_3(U))$ we obtain that \eq{ass} holds
for balls centered at $z_n$ and radius $r_n+d(x'_n, {\cal
Z}_3(U))\asymp r_n$.\\
 Furthermore, by exploiting the Lipschitz regularity of $U$,
 we have
 \neweq{claim1}\int_{\partial
B(x_n,\rho_n)}|\nabla U|^2\leq {3\over \rho_n}
\int_{B(x_n,\rho_n)}|\nabla U|^2\endeq
where $\rho_n=(2L/C)r_n$ and $L$ is the Lipschitz constant in
$\Omega'=\{x\in\Omega:d(x,\partial \Omega)\geq d(x_0,\partial
\Omega)/2\}\supset \{x_n\}$.\\
Hence, in the following let us assume that \eq{ass} holds for a choice of
$x_n\in{\cal Z}_3$ and of radii $\rho_n$  satisfying \eq{claim1}
(we denote $\rho_n$ again by $r_n$).\\
Let $r>0$:
by the monotonicity formula with three phases, there exists $C>0$
(independent of $r$) such that $$\prod_{i=1}^{3}{1\over
r^3}\int_{B(x_n,r)}|\nabla \omega_i(x)|^2 dx \leq C $$ for all
$\omega_1:=u_i$, $\omega_2:=u_j$ and $\omega_3:=\sum_{h\not\in\{
i,j\}} u_h$, such that $x_n\in\partial\supp(u_i)\cap\partial\supp(u_j)\cap\partial\supp(u_l)$ for some
$l\not\in\{ i,j\}$.
Then, due to
assumption  \eq{ass}, we deduce that there exist at most two
components, say $u_1$ and $u_2$, such that
\neweq{con1}
\lim_{n\to\infty} {1\over r_n^2}\int_{B(x_n,r_n)}|\nabla
u_i(x)|^2dx>0,\;i=1,2 \hspace{1cm}{1\over
r_n^2}\int_{B(x_n,r_n)}|\nabla u_i(x)|^2dx\to
0\;\;\;\forall\;i\geq 3.\endeq Let us set $$L_n^2:={1\over
r_n^2}\int_{B(x_n,r_n)}|\nabla U(x)|^2dx\;$$ and
consider the sequence of functions
\neweq{riscalate}{U}_n(x)=\frac{1}{L_nr_n}U(x_n+r_nx)\endeq defined in $x\in B(0,1)$.
Note that, since $U\in {\cal S}^*_{M,k}$ by \eq{emme}, then
$U_n\in {\cal S}^*_{M_n,k}$ where $M_n:=\frac{r_n}{L_n}M$.
Then $\int_{B(0,1)}|\nabla U_n|^2=1$ and, by \eq{claim1},
$\int_{\partial B(0,1)}|\nabla U_n|^2$ is bounded too. By Proposition \ref{poincare.convergenze} Part (c), there
exists $\overline{U}\in H^1(B(0,1))$ such that, up to subsequences, $U_n\rightharpoonup
\overline{U}$ and $\overline{U}\neq 0$. Furthermore, by Part (a), we know that at least two components of $U_n$
does not vanish to the limit: comparing with \eq{con1} we have
$$\overline{U}=(\overline{u}_1,\overline{u}_2,0,\dots,0)$$
with $\overline{u_i}\neq 0$ for
$i=1,2$.
Moreover, since $M_n\to 0$, it holds $\overline{U}\in {\cal S}^*_{0,k}$. This in particular implies that
$\overline{u}_1-\overline{u}_2$ is
harmonic.\\
Now, if $\widehat{\overline{u}}_1(0)>0$, (resp. $\widehat{\overline{u}}_2(0)>0$) then
$\overline{u}_1>0$ (resp. $\overline{u}_2>0$)
 in $B(0,\bar{r})$ for some $\bar{r}>0$. This implies
$\int_{B(0,\bar{r})}\widehat{\overline{u}}_{1,n}(x)\geq \a$ for
some $ \a>0$,  since it converges to the $L^1$--norm of
$\overline{u}_1$ in $B(0,\bar{r})$. We can thus apply Remark
\ref{proprieta_media} to the sequence
$(\widehat{\overline{u}}_{1,n})$ obtaining
$\widehat{\overline{u}}_{1,n}(0)\geq \a/2$. But this is in
contradiction with the fact that, by definition,
$\widehat{\overline{u}}_{1,n}(0)=0$.

Now set $v=\widehat{\overline{u}}_1$ and assume that $v(0)=0$.
Standard results on harmonic function (see \cite{hw}) imply that
$v(r,\theta)\sim r^p\cos p(\theta+\theta_0)$ for some $p\geq 1$.
Thus, by the strong convergence in $H^1$ and a diagonal process,
we can assume that, for $n$ large enough and $i=1,2$, there exists
$m_i>0$ such that $u_{n,i}>m_i$ on a circular sector $A_{n,i}=
\{(\rho,\theta): r<\rho<R,
\alpha_{n,i}<\theta<\beta_{n,i}\}\subset \supp(u_{n,i})$. Here we
can assume, for instance,
$\alpha_{n,1}=\alpha<\beta_{n,1}<\alpha_{n,2}<\beta_{n,2}=\beta$
with $\alpha$ and $\beta$ fixed; note that
$\alpha_{n,2}-\beta_{n,1}\to 0$ as $n\to\infty$. Now, since 0 is a
zero of $U_n$ with multiplicity $m(0)\geq 3$, there exists a third
component, say $u_{n,3}$, and a continuous path $\gamma_n:[0,1]\to
B(0,R)$ such that $\gamma_n(0)=0$, $\gamma_n(1)\in\partial
B(0,R)$, $\gamma_n(t)\in\{\beta_{n,1}<\theta<\alpha_{n,2}\}$,
$u_{n,3}(\gamma_n(t))>0$ for all $t\in (0,1]$. Therefore, if we
set $S=\{ (\rho,\theta): r<\rho<R, \alpha<\theta<\beta\}$ and
denote by ${\omega}_{n,i}$ ($i=1,2$) the  connected component of
$\supp(u_{n,i})$ that contains $A_{n,i}$, then $dist(
S\cap\omega_{n,1},S\cap\omega_{n,2})>0$.\\ Now consider
$$T_n(x)=r_nx+x_n$$ and for $x\in T_n(S)$ define

$$w_n(x)=\left\{\begin{array}{l}u_i(x),\hspace{1cm}x\in T_n(\omega_{n,i})\hspace{1cm}i=1,2\\
                                    -u_i(x),\hspace{1cm}x\in
                                    \supp(u_i)\setminus(T_n(\omega_{n,1})\cup T_n(\omega_{n,2}))\hspace{1cm}i=1,...,k
                                    \end{array}\right. $$
We claim that
\neweq{w.sopra}
-\Delta w_n\geq f(w_n)\hspace{1cm}\mbox{ in } T_n(S)\endeq where
$f(x,s):=f_i(x,s)$ for $x\in \supp(u_i)$.\\
 In order to prove this
assertion, let us drop the dependence on $n$. Then fix any
$\phi>0$, $\phi\in C^1_0(T(S))$: we have to prove
$$\int_{T(S)}\left(\nabla w\nabla \phi-f(w)\phi\right)dx >0.$$ For
easier notation, set $A_1=T(\omega_1)$, $A_2=T(\omega_2)$ and
$B=T(S)\setminus(A_1\cup A_2)$. Now take a partition of the unity
in such a way that $\phi(x)=\phi_1(x)+\phi_2(x)$ where $\phi_i>0$
and $\supp(\phi_i)\subset A_i\cup B$.  For $x\in T(S)$ define the
function $$w_1(x)=\left\{\begin{array}{l}u_1(x),\hspace{1cm}x\in
A_1\cup (B\cap\supp(u_1))\\
                                    -u_i(x),\hspace{1cm}x\in
                                    \supp(u_i) \hspace{1cm} i=2,...,k
                                    \end{array}\right. $$
(Analogous definition for $w_2$.) Since $U\in\cal S$ (namely
$-\Delta \widehat u_1\geq f(\widehat u_1)$), it holds
$$\int_{T(S)}\left(\nabla w_1\nabla \phi_1-f(w_1)\phi_1\right)dx
>0$$ and

$$\int_{T(S)}\left(\nabla w_2\nabla \phi_2-f(w_2)\phi_2\right)dx >0.$$
Summing up the two inequalities we obtain
$$\int_{T(S)}\left(\nabla w\nabla \phi-f(w)\phi\right)dx >
-2\sum_{i=1,2}\int_B\left(\nabla u_i\nabla
\phi_i-f(u_i)\phi_i\right) dx.$$

As we are going to prove, each term of the above  sum is negative, and this finally
 proves assertion \eq{w.sopra}. To this aim let
$i=1$, and introduce  a cut--off function $\eta$  such  that
$\eta=1$ on $B\cap\supp(u_1)$ and $\eta=0$ on $T(S)\setminus B$.
Let $\psi=\eta\phi_1\in C^1_0(T(S))$; then
$$\int_B\left(\nabla u_1\nabla \phi_1-f(u_1)\phi_1\right)dx=
\int_{T(S)}\left(\nabla u_1\nabla \psi-f(u_1)\psi\right)dx.$$ Then
the righthand side is negative since $U\in\cal S$ (namely $-\Delta
u_1\leq f(u_1)$).

{\it Final step.} Now fix $t_n>0$ and $R>r_n>r$ in such a way
that, if we set $y_n=\gamma_n(t_n)$ then it holds $\int_{\partial
B(y_n,r_n)}w_n\geq \a$ for some positive $\a$. (This is due to the
fact that $\partial B(y_n,r_n)\subset
\omega_{n,1}\cup\omega_{n,2}$ except for a small piece of total
length less then $2R(\alpha_{n,2}-\beta_{n,1})\to 0$ as
$n\to\infty$.) This allows to apply Remark \ref{proprieta_media}
to the sequence of rescaled functions $(\frac{w_n}{r_nL_n})$ (as in \eq{riscalate}), since
they satisfy $\Delta \frac{w_n}{r_nL_n}\geq -M_n\to 0$. This gives
$w_n(y_n)>0$, in  contradiction with the fact that, by
construction, $w_n(y_n)=-u_{n,3}(y_n)<0$.
\endproof

\subsection{Local properties of the free boundary}\label{free.boundary}
Let again $N=2$; the purpose of this section is to deepen the investigation
of the behavior of the elements $U$ of the class $\cal S$ around multiple points,
that is, $x\in\Omega$ with $m(x)\geq 2$.
Our main goal is to prove that, near a multiple point, $U$ and its
null set exhibit the same qualitative behavior of  harmonic
functions and their nodal sets. We refer the readers to the
fundamental papers of Alessandrini \cite{alex1} and Hartman Winter
\cite{hw} for the main results about the zero set of harmonic
functions and, more in general, of functions in the kernel of a
divergence type operator.
We shall obtained the desired description for those $U\in\cal S$
with the property that each component
has connected support, i.e.
\neweq{supp.conn} U\in
{\cal S}\mbox{ such that }\supp(u_i) \mbox{ is  connected } \forall
i.\endeq
In order to simplify the topological arguments involved in the proof of our results, in what
follows we shall
make the additional (but not necessary) assumptions that $\O$ is simply connected and
that all the boundary data are not identically zero, namely $\phi_i\not\equiv 0$ on $\partial\O$.
In this situation the assumption \eq{supp.conn} (and Remark \ref{deriv.norm}) implies that
$\supp(u_i)$ is simply connected.\\

Let us start our description by
analyzing points of multiplicity two, that will be denoted by
$${\cal Z}^2=\{x\in \Omega : m(x)=2\}.$$
We already noticed in Remark \ref{deriv.norm} that this set is
locally a regular ${\cal C}^1$ arc around those points where
$\nabla U$ does not vanish. Our first result states that in fact
this is always true:
\begin{lemma}\label{localc1}
Let $x_0\in \Omega$ such that $m(x_0)=2$. Then $\nabla U(x_0)\neq
0$ and ${\cal Z}^2$ is locally a $C^1$--curve through $x_0$.
\end{lemma}
\proof since $x_0\in \partial\{u_1>0\}\cap\partial\{u_2>0\}$, then
for all $r$ small enough  $B(x_0,r)\cap \supp(u_i)=\emptyset$ for
all $i>2$. Assume by contradiction $\nabla u_1(x_0)=0=\nabla
u_2(x_0)$; then $u=u_1-u_2$ is $C^1$ and satisfies the equation
$-\Delta u=a(x) u$ (with
$a(x)=\frac{f_1(u_1(x))-f_2(u_2(x))}{u_1(x)-u_2(x)}$) in
$B(x_0,r)$. By  Lemma \ref{ker.div} there exists a positive,
regular function $\f$ which is radial with respect to $x_0$ such
that
$$-\div\left(\f^2\nabla\left(\frac{u}{\f}\right)\right)= 0$$
on $B(x_0,r)$.
Then $u/\f$ satisfies all the assumptions necessary
to apply the main theorem in \cite{alex1}, which says that the
null level set of $u/\f$ (indeed of $u$) near $x_0$ is made up by a finite number of
curves starting from $x_0$. Obviously in our situation such number
must be even. Now recall that each $\supp(u_i)$ is  connected in
$\Omega$: by a geometrical argument we can see that the null level
set is made up by (two semi--curves joining in) one $C^1$--curve.
But again applying \cite{alex1} we have $\nabla U(x_0)\neq 0$, a
contradiction.
\endproof

\begin{lemma}\label{aderenza} Let $x_0\in{\cal Z}_3$.
Then there exists $\{x_n\}\subset\Omega$ such that $m(x_n)=2$ and
$x_n\to x_0$ \end{lemma} \proof assume not, then there would be an
element $y_0$ of ${\cal Z}_3$ having a positive distance $d$ from
${\cal Z}^2$. Let $r<d/2$: then the ball $B(y_0,r)$ intersects at
least three supports; therefore there exist, say, $x\in\supp(u_i)$
and $z_0\in{\cal Z}_3$ such that $\rho=d(x,z_0)=d(x,{\cal
Z}_3)<d(x,{\cal Z}^2)$. Then the ball $B(x,\rho)$ is tangent from
the interior of $\supp(u_i)$ to ${\cal Z}_3$ in $z_0$; furthermore
$u_i$ solves an elliptic PDE and it is positive on its support: we
thus infer from the Boundary Point
 Lemma that $\nabla u_i(z_0)\neq0$, in contrast with Theorem \ref{grad.nullo}.\endproof

Let us now prove  an asymptotic formula describing the behavior of
$\sum u_i$ in the neighborhood of a multiple point which is
isolated in ${\cal Z}_3$.

\begin{teo}\label{asym} Let $x_0\in \O$ with $m(x_0)=h\geq 3$
that is isolated in ${\cal Z}_3$. Then there exists
$\theta_0\in (-\pi,\pi]$ such that
$$U(r,\theta)=r^{\frac{h}{2}}|\cos(\frac{h}{2}(\theta+\theta_0))|+o(r^{\frac{p}{2}})$$
as $r\to 0$, where $(r,\theta)$ denotes a system of polar coordinates around
$x_0$.
\end{teo}
\proof by assumption  $x_0$ is isolated in ${\cal Z}_3$, then
there is $B=B(x_0,\rho)$ such that ${\cal Z}_3\cap B=\{x_0\}$; furthermore since
each $ \supp(u_i)$ is simply connected, then
$\supp(u_i)\cap\partial B=b_i$ is a connected arc on $\partial B$ for
the $h$ indices involved in $x_0$.
Choosing a slightly smaller radius we can suppose that the
intersection of $\partial B$ with ${\cal Z}^2$ is transversal
(recall that, by Lemma \ref{localc1}, ${\cal Z}^2$ is locally a
$C^1$--curve). Let us assume that $h$ is even:
then we define a function $v(r,\theta)$ such that
$|v(r,\theta)|=|u(r,\theta)|$ and $sign(v(r,\theta))=(-1)^j$ if
$(\rho,\theta)\in b_j$, $j=1,...,h$. Note that the resulting
function is alternately positive and negative on the consecutive
(with respect to $\theta$) local components of $U$. If on the
contrary $h$ is odd, we define $|v(r,\theta)|=|u(r^2,2\theta)|$
and we prescribe an alternating sign to the local components of
$u(\rho^2,2\theta)$. It is worthwhile noticing that the resulting
function $v$ is of class $C^1$ in $\tilde B=B(x_0,\rho^2)$: indeed,
$\supp(u_i)\cap B$ is simply connected for every $i$ and thus
each connected component of $B\setminus u^{-1}(0)$ corresponds to
two components of $\tilde B\setminus v^{-1}(0)$ to which we give
opposite sign. In both the even and the odd cases  $v$ is of class
$C^1$ and it solves an equation of type $-\Delta v=a(x) v$ in
$B\setminus\{x_0\}$ (resp. $\tilde B\setminus\{x_0\}$),  where
$a\in L^\infty$ is given by $f(v)/v$ and $r^2f(v)/v$ respectively.
Moreover $\nabla v(x_0)=0$ by Theorem \ref{grad.nullo}: this
implies that $v$ is in fact solution of the equation on the whole
of $B$ (resp. $\tilde B$) and thus it is of class $C^{2,\a}$. Now,
choosing $r$ small enough, by Lemma \ref{ker.div} we have a
positive, regular function $\f$ such that
$$-\div\left(\f^2\nabla\left(\frac{v}{\f}\right)\right)= 0.$$
Then, we can apply to $v/\f$ the asymptotic formula of Hartman and
Winter as recalled in \cite{hw}. To complete the proof, let us
observe that $h$
represents the number of connected components of $\cup \supp(u_i)$ in a ball centered in $x_0$,
providing the choice of the periodicity of the cosine in the representation formula.
\endproof

The last part of this section is devoted to prove that ${\cal
Z}_3$ consists of (a finite number of) isolated points: thus
Theorem \ref{asym} will provide the complete description of $U$
near multiple points.

We start by proving an intermediate result:

\begin{prop}\label{finite}
${\cal Z}_3$ has a finite number of
connected components.
\end{prop}
\proof
 let $\omega_i:=\supp(u_i)$; take an index pair $(i,j)$ such that $\partial \omega_i$,
$\partial\omega_j$ do intersect and we consider
$$\Gamma_{i,j}=\partial \omega_i\cap\partial \omega_j\cap {\cal Z}^2$$
$$\omega_{i,j}=\omega_i\cup\omega_j\cup\Gamma_{i,j}\;.$$
Since by Lemma \ref{localc1}  $\Gamma_{i,j}$ is locally a regular
arc  and each $\omega_i$ is open, it easily follows that
$\omega_{i,j}$ is open. Furthermore $\omega_{i,j}$ is simply
connected. Indeed, let us consider a loop $\gamma$ in
$\omega_{i,j}$ which is not contractible in $\omega_{i,j}$. This
means that $\Omega_\gamma$, the internal region of the loop, contains at
least a multiple intersection point. Since each support is
connected, there exists $\omega_h$, $h\neq i,j$ such that
$\omega_h\subset \Omega_\gamma$. But this is in contradiction with
the fact that $\omega_h\cap\partial\Omega=\supp(\phi_h)$.\\ A similar
reasoning allows to prove that each $\Gamma_{i,j}$ consists in a
single $C^1$--arc: as straightforward consequence of these facts
the set of multiple points can have only a finite number of
connected components.
\endproof

 Now  we will need  the following definition of adjacent
supports:
\begin{defin} We say that $\omega_i$ and $\omega_j$ are adjacent
if $$\Gamma_{i,j}\neq\emptyset.$$
\end{defin}
Let us list some basic properties:
\begin{itemize}
\item[1.] Every $\omega_i$ is adjacent to some other $\omega_j$.
This follows from the Boundary Point Lemma.
\item[2.] Let us pick $k$ points $x_i\in\omega_i$, $i=1,\dots,k$. If $\omega_i$ and
$\omega_j$ are adjacent, $i<j$, then there exists a smooth arc
$\gamma_{ij}$ with $\gamma_{ij}(0)=x_i$, $\gamma_{ij}(1)=x_j$
lying in $\omega_i\cup\omega_j\cup\{y_{ij}\}$, for some
$y_{ij}\in{\cal Z}^2$.
\item[3.] We can choose the arcs $\gamma_{ij}$ in a manner that they
are mutually disjoint, except for the extreme points.
\end{itemize}
We call $\cal G$ the graph induced by the arcs $\gamma_{ij}$ and
their endpoints.\\

{\bf Construction of an auxiliary function $v$.}
 Aim of this paragraph is to build up by the components of $U$ a $C^1$ function
carrying a sign law which is compatible with the adjacency
relation and solves an elliptic equation. Let us first assume that\\

\centerline{{\it  $\cal G$ has no loops.}}

In this situation the
graph is a disjoint union of a finite number of branches: we
select one of them. Now we define $v$ as follows: if the index $i$
is not involved in the branch, we set $v\equiv 0$ in $\supp(u_i)$.
Next we follow the selected branch of the graph and prescribe a
sign to each vertex by alternating plus and minus. We then define
$v(x)=\pm u(x)$ according to this sign rule. Taking into account
Remark \ref{deriv.norm} and Theorem
\ref{grad.nullo}, with this procedure we obtain a $C^1$--function on $\Omega$.\\

Moreover $v$ solves an elliptic equation as follows: let us define
$f(x,s):=f_i(x,s)$ if $x\in\supp(u_i)$, $i=1,\dots,k$. Then by
construction
\neweq{equaz.v}-\Delta v=a(x)v\hspace{1cm}\mbox{ in }\Omega\setminus{\cal Z}_3
\endeq
where  $a=f(v)/v\in L^\infty$.  In fact we are going to prove that
$v$ solves the elliptic equation on the whole of $\Omega$. We need
a technical lemma

\begin{lemma} \label{normal}
Let ${\cal Z}_3^\varepsilon=\{x\in\Omega\;:\;m(x)=2\}\cap
B_\varepsilon ({\cal Z}_3)$. For every $\delta>0$ there exists
$\eps>0$ such that $$\int_{{\cal Z}_3^\varepsilon}\left|\nabla
v\right|<\delta.$$
\end{lemma}
\proof by testing the equation \eq{equaz.v} with the test function
$\f=1$ and integrating over the set $u_i>\alpha$ we obtain the
bound, independent of $\alpha$ and $i$,
$$\int_{\partial\{u_i>\alpha\}}\left|\nabla v\right|<C$$ and
therefore, passing to the limit as $\alpha\to 0$,
$$\int_{\partial\{u_i>0\}}\left|\nabla v\right|<C.$$ The assertion
then follows from Lemma \ref{aderenza}, together with Proposition
\ref{finite}.
 \endproof

\begin{lemma}\label{proalex} $v$ solves $-\Delta v=a(x)v$ in $\Omega$.
\end{lemma}
\proof let us fix a connected component of ${\cal Z}_3$, named $A$.
Thanks to Lemma \ref{aderenza}, for any $\eps$ we can take
a neighborhood of $A$, say ${\cal V}_\eps\subset B_\varepsilon ({\cal Z}_3)$,
 in such a way that the boundary $\partial{\cal
V}_\eps$ is the union of a finite number of arcs of
${\cal Z}^2$ and supplementary union of pieces of
total length smaller than $C\eps$. Let $\varphi$ be a test
function. We write

$$\begin{array}{c}
 \int_\Omega(\nabla v\nabla\varphi-a(x)v\varphi)=\\
 =\int_{\Omega\setminus{\cal V}_\eps}(\nabla v\nabla\varphi-
 a(x)v\varphi)+\int_{{\cal V}_\eps}(\nabla
 v\nabla\varphi-a(x)v\varphi)\leq\\
 \leq C\int_{\partial{\cal V}_\eps}|\nabla v|+C\int_{{\cal V}_\eps}(|\nabla
 v|+|v|).
 \end{array}$$

Let $\delta>0$: we can find $\eps>0$ such that Lemma \ref{normal}
holds. Moreover, by Theorem \ref{grad.nullo}, we can assume that $\eps$
is taken so small that $\sup_{{\cal V}_\eps}(|\nabla
v|+|v|)<\delta$. Hence the above integral is bounded by $C\delta$.
Since $\delta$ was arbitrarily chosen we obtain that $v$ solves
the equation in a distributional sense. Usual regularity arguments
allow us to complete the proof.
 \endproof

With this  we easily deduce the final result, that is, $U$ has
only a finite number of multiple points:
\begin{lemma} \label{interior} The set ${\cal Z}_3$ consists of a
finite number of points.
\end{lemma}
\proof since  by Lemma \ref{proalex} $v$ is  a solution of
$-\Delta v=a(x)v$ in $\Omega$, locally we can reduce $v$ to a
function in the kernel of a divergence--type operator as in Lemma
\ref{ker.div}. Then the results of \cite{alex1} and \cite{hw}
ensures that  $v$ has only a finite number of multiple points.
\endproof

Let us now go back to the case when the graph associated to $U$ presents loops:\\

\centerline{{\it  $\cal G$ has a loop.}}

Let us first define an order relation between loops, according
whether one is contained in the interior region of the other. Let
us select a {\it minimal} loop  $\gamma$ (no other loops are
contained in $\Omega_\gamma$, its interior region). We can assume
that $\Omega_\gamma$ contains at least an element $x_0\in{\cal
Z}_3$ (if not, we can perform a conformal inversion exchanging the
inner with the outer points): fix $x_0$ as the origin. Note that
the supports  involved in a minimal loop have the remarkable
property that, for all $i$, $\omega_i$ is adjacent to $\omega_j$
for only two indices different from $i$. Thanks to this property,
if the number of vertex of $\gamma$ is even, we manage to assign a
sign law to all the subset of $\cal G$ contained in
$\Omega_\gamma$ so that adjacent supports have opposite sign: it
suffices to follow the loop and prescribe alternating sign to its
vertices. We
define $v(x)=\pm u(x)$ according to this law.\\
If the number of vertex of $\gamma$ is odd, we wish to ``double''
the loop. To this aim, we define new $\omega_i$'s by taking the
complex square roots of the old ones. In this way the new loop
$\gamma$ will have an even number of edges and we define
$v(r,\theta)=\pm u(r^2,2\theta)$ by giving  alternating sign at
the vertices of the loop.
In this way we define a function $v$ in $\Omega_\gamma$ which is of class
$C^1$ thanks to Remark \ref{deriv.norm} and Theorem
\ref{grad.nullo}. Again, $v$  solves the elliptic equation
\neweq{equaz.v.2}-\Delta v=a(x)v\hspace{1cm}\mbox{ in }\Omega_\gamma\setminus{\cal Z}_3
\endeq
where $a$ is either $a=f(v)/v$ or $a=r^2f(v)/v$
according to the construction of $v$.
As in the previous case, it is possible to show that $v$ solves the
equation on the whole of $\Omega_\gamma$. The proof follows exactly that of Lemma \ref{proalex},
with two remarks; first note that Lemma \ref{normal} still holds. Furthermore, in this situation the set
 ${\cal Z}_3\cap \Omega_\gamma$ is
connected (assume not; then in $\O_\gamma$ there are two points of ${\cal Z}_3$
which are connected by a regular arc of double points. By a simple geometrical argument
in the plane this implies that one of the supports is adjacent to other three different
supports: as already observed, this is in contradiction with the minimality of the loop).
Then ${\cal V}_\eps$ will be a neighborhood of the whole ${\cal Z}_3\cap \Omega_\gamma$ with the properties
required in the proof of that lemma.

This immediately provides
\begin{lemma} \label{interior.2} There is only one point of ${\cal Z}_3$ lying in the interior
of a minimal loop $\gamma$.
\end{lemma}
\proof let $\Omega_\gamma$ denote the internal region of $\gamma$:
since  by Lemma \ref{proalex} $v$ is  a solution of
$-\Delta v=a(x)v$ in $\Omega\gamma$,  locally we can reduce $v$ to a
function in the kernel of a divergence--type operator as in Lemma \ref{ker.div}.
Then by \cite{alex1} and \cite{hw} we know that only a
finite number of points of ${\cal Z}_3$ lie in $\Omega_\gamma$.
On the other hand, by the minimality of $\gamma$, ${\cal
Z}_3\cap\Omega_\gamma$ is connected.
Thus the origin is the only one multiple point contained in
$\Omega_\gamma$.
\endproof

We can finally prove that  multiple points are
isolated
\begin{teo} \label{isolz3.2} The set ${\cal Z}_3$ consists of
isolated points.
\end{teo}
\proof recalling Proposition \ref{finite}, we argue by induction
over the number $h$ of connected components of the set ${\cal
Z}_3$. If $h=1$ then, by Lemma \ref{interior.2} there is at most one
minimal loop of the adjacency relation. If there is one, then
Lemma \ref{interior.2} gives the desired assertion. If there are
none,  the thesis directly follows by Lemma \ref{interior}. Now, let the Theorem be true for $h$
and assume that ${\cal Z}_3$ has $h+1$ connected components.
Again, if the adjacency relation has no loops we are done.
Otherwise, we apply Lemma \ref{interior.2} to treat those connected
components contained in the interior of the minimal loop and the
inductive hypothesis to treat all those contained in the outer
region.\endproof
\begin{rem}{\rm
Having proved that the multiple points are isolated, the existence of points of multiplicity
zero can be easily ruled out for connected domains. }
\end{rem}

\end{document}